\documentclass[12pt, reqno]{amsart}
\usepackage{latexsym}
\usepackage{mathtools}
\usepackage{amsmath}
\usepackage{amsthm}
\usepackage{mathrsfs}
\usepackage{lmodern}
\usepackage[T1]{fontenc}
\usepackage{textcomp}
\usepackage[utf8]{inputenc}
\usepackage{color}
\usepackage{amssymb}
\usepackage{enumerate}
\usepackage[abbrev]{amsrefs}
\usepackage{indentfirst}
\usepackage{graphicx}
\usepackage{bmpsize}
\usepackage{hyperref}
\usepackage{enumerate}
\usepackage{enumitem}
\usepackage{url}
\usepackage{autobreak}
\usepackage[justification=justified]{caption}
\usepackage{bm}
\usepackage{here}
\renewcommand{\thefootnote}{} 

\theoremstyle{plain} 
\newtheorem{theorem}{\indent\sc Theorem}[section]
\newtheorem{lemma}[theorem]{\indent\sc Lemma}
\newtheorem{corollary}[theorem]{\indent\sc Corollary}
\newtheorem{proposition}[theorem]{\indent\sc Proposition}

\theoremstyle{definition} 
\newtheorem{definition}[theorem]{\indent\sc Definition}
\newtheorem{remark}[theorem]{\indent\sc Remark}
\newtheorem{example}[theorem]{\indent\sc Example}

\newenvironment{sproof}{%
  \proof}{\endproof}
\makeatletter
  
  \@addtoreset{equation}{section}
\makeatother
\newcommand{\N}{P_n}
\newcommand{\NN}{\mathfrak{C}_n}
\newcommand{\Nn}{\N^{<\omega}}
\newcommand{\Pset}{\bm{P}}
\newcommand{\Qset}{\bm{Q}}
\newcommand{\Pn}{\Pset^{<\omega}}
\newcommand{\Qn}{\Qset^{<\omega}}
\newcommand{\x}[2]{{x_{#1}}_{[#2]}}
\newcommand{\tx}[2]{{{\tilde{x}}_{#1}}{}_{[#2]}}
\newcommand{\hx}[2]{{\hat{x}_{#1}}{}_{[#2]}}
\newcommand{\X}[2]{X_{#1}^{#2}}

\newcommand{\y}[2]{y_{#1}^{#2}}
\newcommand{\restr}[1]{|_{#1}}
\newcommand{\supp}[1]{\textrm{supp}(#1)}
\begin{document}

\keywords{amenable, finitely presented, Thompson group, torsion free}

\title[An $n$-adic generalization of the Lodha--Moore group]
{AN $n$-ADIC GENERALIZATION OF THE LODHA--MOORE GROUP}
\author{YUYA KODAMA}
\date{}
\renewcommand{\thefootnote}{\arabic{footnote}}  
\setcounter{footnote}{0} 
\begin{abstract}
We generalized the Lodha--Moore group $G_0$, similar to the case of the Thompson group $F$, into an $n$-adic version $G_0(n)$. 
This generalization is given by regarding the Lodha--Moore group as the group consisting of homeomorphisms on the Cantor space. 
Just as the case of $F$ and $F(n)$, we also showed that $G_0(n)$ and $G_0$ share many properties. 
\end{abstract} 
\maketitle
\section{Introduction}
In \cite{lodha2016nonamenable}, Lodha and Moore introduced the group $G_0$ consisting of piecewise projective homeomorphisms of the real projective line. 
This group is a finitely presented torsion free counterexample to the von Neumann-Day problem \cite{Neumann1929, day1950means}, which asks whether every nonamenable group contains nonabelian free subgroups. 
Although counterexamples to the von Neumann-Day problem are known \cite{MR586204, adian1979burnside, MR682486, ol2003non, ivanov2005embedding,  monod2013groups}, it is still an open question whether the Thompson group $F$ can be a new counterexample. 

The Lodha--Moore group $G_0$ has similar properties to the Thompson group $F$. 
Indeed, the Lodha--Moore group can also be defined as a finitely generated group consisting of homeomorphisms of the space of infinite binary sequences, whose generating set is obtained by adding an element to the well-known finite generating set of $F$. 
Both have (small) finite presentations \cite{brown1987finiteness, lodha2016nonamenable}, normal forms with infinite presentations \cite{brown1984infinite, lodha2020nonamenable}, simple commutator subgroups \cite{cannon1996introductory, burillo2018commutators}, trivial homotopy groups at infinity \cite{brown1984infinite, zaremsky2016hnn}, no nonabelian free subgroups \cite{brin1985groups, lodha2016nonamenable, monod2013groups}, and are of type $F_\infty$ \cite{brown1984infinite, lodha2020nonamenable}. 

On the other hand, there exist various generalizations of the Thompson group $F$. 
One of the most natural ones is the $n$-adic Thompson group $F(n)$, which is obtained by replacing infinite binary sequences with infinite $n$-ary sequences. 
Even now, it is still being actively studied whether what is true for the group $F$ is also true for the generalized group and what interesting properties can be obtained under the generalization \cite{MR3609402, MR3565428, MR3741884, sheng2021divergence, MR4448422}. 

In this paper, we generalize the Lodha--Moore group similarly and study its properties. 
Namely, we define an $n$-adic generalized group $G_0(n)$ of the Lodha--Moore group $G_0$ and show that several properties which hold for $G_0$ also hold for $G_0(n)$. 
We remark that $G_0(2)$ is isomorphic to $G_0$. 

Let $n, m \geq 2$. 
We show the following: 
\begin{theorem}
\begin{enumerate}[font=\normalfont]
\item The group $G_0(n)$ admits an infinite presentation with a normal form of elements. 
\item The group $G_0(n)$ is finitely presented. 
\item The group $G_0(n)$ is nonamenable. 
\item The group $G_0(n)$ has no free subgroups. 
\item The group $G_0(n)$ is torsion free. 
\item The groups $G_0(n)$ and $G_0(m)$ are isomorphic if and only if $n=m$ holds. 
\item The commutator subgroup of the group $G_0(n)$ is simple. 
\item The center of the group $G_0(n)$ is trivial. 
\item There does not exist any nontrivial direct product decomposition of the group $G_0(n)$.  
\item There does not exist any nontrivial free product decomposition of the group $G_0(n)$. 
\end{enumerate}
\end{theorem}

This paper is organized as follows. 
In Section \ref{section_F(n)}, we recall the definition and properties of the generalized Thompson $F(n)$. 
In Section \ref{section_G0(n)}, we first recall the definition of the group $G_0$ and define the group $G_0(n)$. 
Then by using presentations of $F(n)$ given in Section \ref{section_F(n)}, we define a normal form of elements in $G_0(n)$ and give an infinite presentation of $G_0(n)$. 
Finally, in Section \ref{section_G0(n)_properties}, we study several properties of $G_0(n)$. 

Let us mention some open problems which are known to hold in the case of $G_0$. 
First, it is an interesting question whether $G_0(n)$ can be realized as a subgroup of the group of piecewise projective homeomorphisms of the real projective line (Monod's group $H$ \cite{monod2013groups}). 
The second problem is whether this group is of type {$\rm F_\infty$}, and all homotopy groups are trivial at infinity.  
If it has these two properties, then ${G_0(n)}$ is an example of an (infinite) family of groups satisfying all Geoghegan's conjectures for the Thompson group $F$. 

Furthermore, we can consider some groups related to $G_0(n)$. 
The first one is constructed by using another definition of the map $y$ defined in Section \ref{subsection_def_G_0(n)}. 
Although we define the map so that $G_0$ is naturally a subgroup of $G_0(n)$, we can consider several different generalizations. 
In particular, there may exist a subgroup of Monod's group $H$, which should be called a ``generalized Lodha--Moore group'' rather than our group. 

We can also construct groups that contain $G_0(n)$. 
In \cite{lodha2016nonamenable}, the group $G$ is defined, where $G$ contains $G_0$ as a subgroup. 
For our group, by adding some of the generators $y_0$, $y_{(n-1)1}, \dots, y_{(n-1)(n-2)}, y_{(n-1)}$ where each $(n-1)i$ is a concatenation of $(n-1)$ and $i$, we can define not only the group $G(n)$, which corresponds to $G$ but also the groups ``between'' $G_0(n)$ and $G(n)$. 
\section{The generalized Thompson group $F(n)$}\label{section_F(n)}
\subsection{Definition} \label{subsection_F(n)_definition}
Let $n \geq 2$. 
There exist several ways to define the Brown--Thompson group $F(n)$ \cite{brin1985groups, brown1987finiteness, burillo2001metrics, guba1997diagram}. 
In this paper, we define it as a group of homeomorphisms on the $n$-adic Cantor set. 
We use tree diagrams to represent elements of the group visually. 

We define $\N$ to be the set $\{0, 1, \dots, n-1 \}$. 
We endow $\N$ with the discrete topology and endow $\NN=P_n^\omega=P_n \times P_n \times \cdots$ with the product topology. 
Note that $\NN$ and the Cantor set are homeomorphic. 
We also consider the set of all finite sequences on $\N$ and write $\Nn$ for it. 
For $s \in \Nn$ and $t \in \Nn$ (or $\NN$), the concatenation is denoted by $st$. 
The group $F(n)$ is a group generated by the following $n$  homeomorphisms: 
\begin{align*}
x_0&\colon \NN \to\NN;
\begin{cases}
0k\eta \mapsto k\eta &(k<n-1) \\
0(n-1)\eta \mapsto (n-1)0\eta \\
k\eta \mapsto (n-1)k \eta &(0<k<n),
\end{cases} \\
x_1&\colon \NN \to \NN;
\begin{cases}
0\eta \mapsto 0\eta \\
1k \eta \mapsto k\eta &(k<n-2) \\
1(n-2)\eta \mapsto (n-1)0\eta \\
1(n-1)\eta \mapsto (n-1)1 \eta \\
k\eta \mapsto (n-1)k\eta &(1<k<n),
\end{cases} \\
&\vdots \\
x_{n-2}&\colon\NN \to \NN;
\begin{cases}
k \eta \mapsto k \eta &(k<n-2) \\
(n-2)0\eta \mapsto (n-2)\eta \\
(n-2)k \eta \mapsto (n-1)(k-1) &(0<k\leq n-1) \\
(n-1)\eta \mapsto (n-1)(n-1) \eta
\end{cases}
\shortintertext{and}
\x{0}{(n-1)}&\colon \NN \to \NN;
\begin{cases}
k\eta \mapsto k \eta &(k<n-1) \\
(n-1)\eta \mapsto (n-1)x_0(\eta) 
\end{cases}
\end{align*}

\begin{figure}[tbp]
	\centering
	\includegraphics[width=150mm]{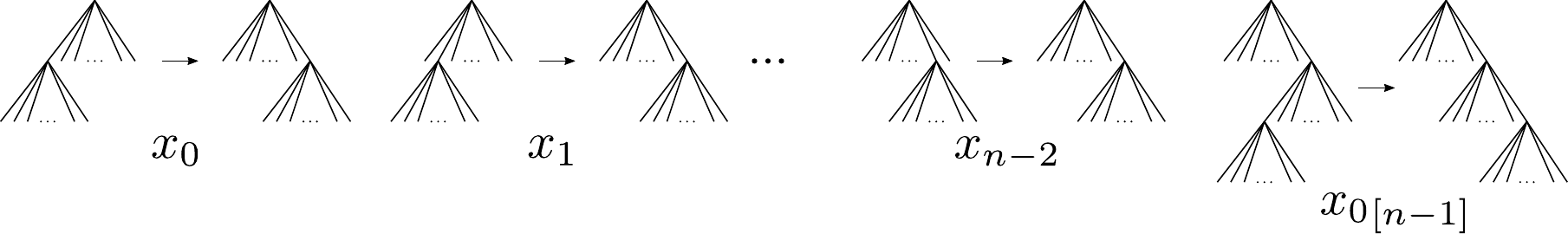}
	\caption{Tree diagrams of the homeomorphisms. }
	\label{generator_Fn}
\end{figure}

\begin{figure}[tbp]
	\centering
	\includegraphics[width=30mm]{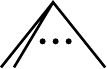}
	\caption{An $n$-caret. }
	\label{n-caret}
\end{figure}

These maps are represented by tree diagrams as in Figure \ref{generator_Fn}. 
Here, we briefly review the definition of tree diagrams. 
See \cite{burillo2001metrics} for details. 
A \textit{rooted $n$-ary tree} is a finite tree with a top vertex (\textit{root}) with $n$ edges, and all vertices except the root have degree only $1$ (\textit{leaves}) or $n+1$. 
In the rest of this paper, we call this $n$-ary tree. 
We define a \textit{$n$-caret} to be an $n$-ary tree with no vertices whose degree is $n+1$ (see Figure \ref{n-caret}). 
Then, each $n$-ary tree is obtained by first attaching an $n$-caret to a single vertex, and then iteratively attaching an $n$-caret to a certain leaf of the obtained $n$-ary tree. 
We always assume that the root is the top and the others are descendants. 

Each $n$-ary tree can be regarded as a finite subset of $\Nn$. 
To do this, we label each edge of each $n$-caret by $0, 1, \dots, n-1$ from the left. 
Since every leaf corresponds to a unique path from the root to the leaf, we can regard it as an element in $\Nn$. 

Let $T_+$ and $T_-$ be $n$-ary trees with $m$ leaves. 
Let $a_1, \dots, a_m$ be elements in $\Nn$ with lexicographic order corresponding to the leaves of $T_+$. 
For $T_-$, define $b_1, \dots, b_m$ in the same way. 
Then, for every $\zeta \in \NN$, there exists $i$ uniquely such that $\zeta=a_i\eta$ for some $\eta \in \NN$. 
Thus we obtain a homeomorphism $a_i \eta \mapsto b_i \eta$. 
It is known that every homeomorphism obtained from two $n$-ary trees with the same number of leaves in this way is generated by a composition of $x_0, x_1, \dots, x_{n-2}, \x{0}{(n-1)}$. 
See \cite[Corollary 10.9]{meier2008groups} for the case $n=2$. 

Let $i$ in $\{0, \dots, n-2 \}$ and $\alpha$ in $\Nn$. 
We define the map $\x{i}{\alpha}\colon \NN \to \NN$ by 
\begin{align*}
\x{i}{\alpha}(\zeta)&=
\begin{cases}
\alpha x_i(\eta) & (\zeta=\alpha \eta) \\
\zeta & (\zeta \neq \alpha \eta), 
\end{cases}
\end{align*}
and we define 
\begin{align*}
X(n)&\coloneqq
\left\{ \x{i}{\alpha} \mid i = 0, \dots, n-2, \alpha \in \Nn \right\}. 
\end{align*}
This set contains the well-known infinite generating set of $F(n)$. 

Let $X^\prime(n)\coloneqq \{ \x{0}{s}, \dots, \x{n-2}{s} \mid s=(n-1)^i \in \Nn, i\geq0 \}$. 
We denote each element as $x_0=X_0, \dots, x_{n-2}=X_{n-2}, \x{0}{(n-1)}=X_{n-1}, \x{1}{(n-1)}=X_{n}, \dots$ only in this section for the sake of simplicity. 
In this notation, $X^\prime(n)=\{X_0, X_1, \dots\}$ holds. 
For $i<j$, we have $X_i^{-1}X_jX_i=X_{j+n-1}$. 
This implies that every element in $F(n)$ has the following form: 
\begin{align*}
\X{i_1}{r_1}\X{i_2}{r_2}\cdots \X{i_m}{r_m}\X{j_k}{-s_k}\cdots\X{j_2}{-s_2}\X{j_1}{-s_1}
\end{align*}
where $i_1<i_2<\cdots<i_m\neq j_k>\cdots>j_2>j_1$ and $r_1, \dots, r_m, s_1, \dots, s_k>0$. 
We require that this form satisfies the following additional condition: 
if $X_i$ and $X_i^{-1}$ both occur in the form, then one of 
\begin{align*}
X_{i+1}, \X{i+1}{-1}, X_{i+2}, \X{i+2}{-1}, \dots, X_{i+n-1}, \X{i+n-1}{-1}
\end{align*}
also occur in the form. 
It is known that this form with the additional condition always uniquely exists. 
Thus we call this form the \textit{normal form} of elements of $F(n)$. 
The proof for the case $n=2$ is in \cite[Section 1]{brown1984infinite}. 

There exist two well-known presentations for $F(n)$ \cite{brin1998automorphisms, guba1997diagram}: 
\begin{align*}
F(n) &\cong \langle X_0, X_1, X_2, \dots  \mid \text{$X_i^{-1} X_j X_i=X_{j+n-1}$, for $i<j$}  \rangle \\
&\cong \left\langle X_0, X_1, \dots, X_{n-1} \;\middle|\;
\begin{array}{l} \text{${X_k}^{X_0}={X_k}^{X_i}$ ($1 \leq i<k\leq n-1$)}, \\
\text{${X_k}^{X_0X_0}={X_k}^{X_0X_i}$ ($1 \leq i, k\leq n-1$ and $k-1\leq i$)}, \\
{X_k}^{X_0X_0X_0}={X_1}^{X_0 X_0 X_{n-1}}
\end{array}
\right\rangle, 
\end{align*}
where $x^y$ denotes $y^{-1}x y$. 

For any $x \in X(n)$, let $w(x)$ be its normal form. 
We define $R_{F(n)}$ to be the set
\begin{align*}
\{ X_i^{-1} X_j X_i=X_{j+n-1}, x=w(x) \mid 0 \leq i <j, x \in X(n) \}. 
\end{align*}
Then $\langle X(n) \mid R_{F(n)} \rangle$ is isomorphic to $F(n)$. 
\section{The Lodha--Moore group and its generalization}\label{section_G0(n)}
\subsection{The Lodha--Moore group}
In this section, we briefly review the original Lodha--Moore group $G_0$. 
Let $x_0$, $\x{0}{1}$ be the maps defined in Section \ref{subsection_F(n)_definition} for $n=2$. 
The group $G_0$ is generated by these two maps and one more generator called $y_{10}$. 
To define this generator, we first define the homeomorphism called $y$. 
\begin{definition}\label{definiton_y_2}
The map $y$ and its inverse $y^{-1}$ is defined recursively based on the following rule: 
\begin{align*}
&y\colon \mathfrak{C}_2 \to \mathfrak{C}_2 & &y^{-1}\colon \mathfrak{C}_2 \to \mathfrak{C}_2 \\
&y(00\zeta)=0y(\zeta) & &y^{-1}(0\zeta)=00y^{-1}(\zeta) \\
&y(01\zeta)=10y^{-1}(\zeta) & &y^{-1}(10\zeta)=01y(\zeta) \\
&y(1\zeta)=11y(\zeta), & &y^{-1}(11\zeta)=1y^{-1}(\zeta). 
\end{align*}
\end{definition}
For each $s$ in $2^{<\omega}$, we also define the map $y_s$ by setting
\begin{align*}
y_s(\xi)&=
 \left \{
\begin{array}{cc}
s y(\eta), & \xi=s\eta \\
\xi, & \mbox{otherwise}. 
\end{array}
\right. 
\end{align*}
We give an example of a calculation of the $y_{001}$ on $00101101 \cdots$: 
\begin{align*}
y_{001}(00101101\cdots)&=001y(01101\cdots)\\
&=00110y^{-1}(101\cdots)\\
&=0011001y(1\cdots). 
\end{align*}
\begin{definition}
The group $G_0$ is a group generated by $x_0$, $\x{0}{1}$, and $y_{10}$. 
\end{definition}
The group $G_0$ is also realized as a group of piecewise projective homeomorphisms. 
\begin{proposition}[{\cite[Proposition 3.1]{lodha2016nonamenable}}]\label{proposition_piecewiseprojective}
The group $G_0$ is isomorphic to the group generated by the following three maps of $\mathbb{R}$: 
\begin{align*}
a(t)&=t+1, 
&
b(t) &= \left \{
\begin{array}{cc}
t & \mbox{\rm{if} $t \leq 0$} \\
\frac{t}{1-t} & \mbox{\rm{if} $0 \leq t \leq \frac{1}{2}$} \\
3-\frac{1}{t} & \mbox{\rm{if} $\frac{1}{2} \leq t \leq 1$} \\
t+1 & \mbox{\rm{if} $1 \leq t$}, 
\end{array}
\right. 
& 
c(t) &= \left \{
\begin{array}{cc}
\frac{2t}{1+t} & \mbox{\rm{if} $0 \leq t \leq 1$} \\
t & \mbox{\rm{otherwise}}. 
\end{array}
\right. 
\end{align*}
\end{proposition}

This proposition is shown by identifying $\mathfrak{C}_2$ with $\mathbb{R}$ by the following maps: 
\begin{align*}
&\varphi\colon \mathfrak{C}_2 \to [0, \infty] & 
&\Phi\colon \mathfrak{C}_2 \to \mathbb{R} \cup \{ \infty \} \\
&\varphi(0\xi)=\frac{1}{1+\frac{1}{\varphi(\xi)}} & 
&\Phi(0\xi)=-\varphi(\tilde{\xi}) \\
&\varphi(1\xi)=1+\varphi(\xi) & 
&\Phi(1\xi)=\varphi(\xi), 
\end{align*}
where $\mathbb{R} \cup \{ \infty \}$ denotes the real projective line and $\tilde{\xi}$ denotes the element of $\mathfrak{C}_2$ obtained from $\xi$ by replacing all $0$s with $1$s and vice versa. 
For instance, if $\xi=0100\cdots$, then $\tilde{\xi}=1011\cdots$. 
We note that for every $x \in \mathbb{R} \cup \{ \infty \}$, the inverse image $\Phi^{-1}(\{x\})$ is either a one-point set or a two-point set. 
Also note that $G_0$ is a subgroup of Monod's group $H$ \cite{monod2013groups}. 
\subsection{$n$-adic Lodha--Moore group} \label{subsection_def_G_0(n)}
In order to define a generalization of the Lodha--Moore group, we first define a homeomorphism on $\NN={\{0, \dots, n-1 \}}^\omega$ corresponding to $y$ for the case of $n=2$. 
We fix $n \geq 2$ and also denote this map by $y$ as in the case of $n=2$. 
\begin{definition}
The map $y$ and its inverse map $y^{-1}$ is defined recursively based on the following rule: 
\begin{align*}
y\colon \NN &\to \NN & y^{-1}\colon \NN &\to \NN \\
y(00\zeta)&=0y(\zeta)& y^{-1}(0\zeta)&=00y^{-1}(\zeta) \\
y(0k\zeta)&=k\zeta & y^{-1}(k\zeta)&=0k\zeta \\
y(0(n-1)\zeta)&=(n-1)0y^{-1}(\zeta) & y^{-1}((n-1)0\zeta)&=0(n-1)y(\zeta) \\
y(k\zeta)&=(n-1)k\zeta & y^{-1}((n-1)k\zeta)&=k\zeta \\
y((n-1)\zeta)&=(n-1)(n-1)y(\zeta) & y^{-1}((n-1)(n-1)\zeta)&=(n-1)y^{-1}(\zeta)
\end{align*}
where $k$ is in $\{1, \dots, n-2\}$. 
\end{definition}
\begin{remark}
Note that if we restrict the domain to $\{0, n-1\}^\omega$, $y$ is exactly the map defined in Definition \ref{definiton_y_2} under the identification of $n-1$ with $1$. 
We will use this fact to reduce the discussion to the case of $n=2$. 
\end{remark}
For each $s$ in $\Nn$, define the map $y_s$ by setting
\begin{align} \label{n-adicy_definition}
y_s(\xi)&=
 \left \{
\begin{array}{cc}
s y(\eta), & \xi=s\eta \\
\xi, & \mbox{otherwise}. 
\end{array}
\right. 
\end{align}
\begin{definition}
We define $G_0(n)$ to be the group generated by the $n+1$ elements $x_0, \dots, x_{n-2}, {x_0}_{[n-1]}$, and $y_{(n-1)0}$. 
We call this group the \textit{$n$-adic Lodha--Moore group}. 
\end{definition}
For an infinite presentation of $G_0(n)$ described in Corollary \ref{G_0(n)_presentation}, we introduce an infinite generating set of this group. 
Let $i$ in $\{0, \dots, n-2 \}$ and $\alpha$ in $\Nn$. 
We recall that the map $\x{i}{\alpha}$ is defined as follows: 
\begin{align*}
\x{i}{\alpha}(\zeta)&=
\begin{cases}
\alpha x_i(\eta) & (\zeta=\alpha \eta) \\
\zeta & (\zeta \neq \alpha \eta). 
\end{cases}
\end{align*}
We also recall that the group $F(n)$ is generated by the following infinite set: 
\begin{align*}
X(n)=
\left\{ \x{i}{\alpha} \mid i = 0, \dots, n-2, \alpha \in \Nn \right\}. 
\end{align*}
In addition, let 
\begin{align*}
Y(n)\coloneqq
\left\{ y_\alpha \;\middle|\;
\begin{array}{l} \alpha \in \Nn, \\
\mbox{$\alpha \neq \epsilon, 0^i, (n-1)^i $ for any $i \geq1$} \\
\mbox{the sum of each number in  $\alpha$ is equal to $0 \bmod {n-1}$} 
\end{array}
\right\}, 
\end{align*}
where $\epsilon$ denotes the empty word. 
We remark that for $\alpha=\alpha_1 \cdots \alpha_m \in \Nn$, the natural actions of $x_0, \dots, x_{n-2}, \x{0}{n-1}$, and $y$ preserve the value of $\alpha_1+\cdots+\alpha_m \bmod {n-1}$. 
The set $Z(n)\coloneqq X(n) \cup Y(n)$ also generates the group $G_0(n)$. 

\subsection{Infinite presentation and normal form}\label{infinite_presentation}
In this section, we give a unique word with ``good properties'' for each element of $G_0(n)$ (Definition \ref{normal_form}). 
In this process, we also give an infinite presentation of $G_0(n)$ (Corollary \ref{G_0(n)_presentation}). 
Almost all results in the rest of this section follow along the lines of the arguments in \cite{lodha2020nonamenable, lodha2016nonamenable}. 
Note that if we only want to get an infinite presentation of $G_0(n)$, it is shorter to just follow \cite{lodha2016nonamenable}. 

Let $s $ in $\Nn$ and $t$ in $\Nn$ or $\NN$. 
We write $s \subset t$ if $s$ is a proper prefix of $t$ and write $s \subseteq t$ if $s \subset t$ or $s=t$. 
We say that $s$ and $t$ are independent if one of the following holds: 
\begin{itemize}
\item $s, t \in \Nn$ and neither $s \subseteq t$ nor $t \subseteq s$ holds. 
\item $s \in \Nn$, $t \in \NN$ and $s$ is not any prefix of $t$. 
\item $s, t \in \NN$ and $s\neq t$. 
\end{itemize}
In all cases, we write $s \perp t$. 

Let $s(i), t(i)$ denote the $i$-th number of $s, t$, respectively. 
Then we say $s<t$ if one of the following is true: 
\begin{itemize}
\item[(a)] $t \subset s$; 
\item[(b)] $s \perp t$ and $s(i)<t(i)$, where $i$ is the smallest integer such that $s(i) \neq t(i)$. 
\end{itemize}
We note that this order is transitive. 
For elements in $\NN$, we use the same symbol $<$ to denote the lexicographical order. 

We claim that the following collection of relations gives a presentation of $G_0(n)$ (Corollary \ref{G_0(n)_presentation}): 
\begin{enumerate}
\item the relations of $F(n)$ in $R_{F(n)}$; 
\item $y_t\x{i}{s}=\x{i}{s}y_{\x{i}{s}(t)}$ for all $i$ and $s, t \in \Nn$ such that $y_t \in Y(n)$, and $\x{i}{s}(t)$ is defined; 
\item $y_s y_t =y_t y_s$ for all $s, t \in \Nn$ such that $y_s, y_t \in Y(n)$, and $s\perp t$; 
\item $y_s=\x{0}{s} y_{s0} y_{s(n-1)0}^{-1}y_{s(n-1)(n-1)}$ for all $s \in \Nn$ such that $y_s \in Y(n)$. 
\end{enumerate}
We note that $\x{0}{n-1}((n-1)0)$ is not defined, for example. 
All relations can be verified directly. 
We write $R(n)$ for the collection of these relations. 

First, we define a form that makes it easy to compute the composition of maps. 
\begin{definition}
A word $\Omega$ on $Z(n)$ is in \textit{standard form} if $\Omega$ is a word such as $f\y{s_1}{t_1}\cdots\y{s_m}{t_m}$ where $f$ is a word on $X(n)$ and $\y{s_1}{t_1}\cdots\y{s_m}{t_m}$ is a word on $Y(n)$ with the condition that $s_i<s_j$ if $i<j$. 
\end{definition}
In some cases, it is helpful to use the following weaker form. 
\begin{definition}
A word $\Omega$ on $Z(n)$ is in \textit{weak standard form} if $\Omega$ is a word such as $f\y{s_1}{t_1}\cdots\y{s_m}{t_m}$ where $f$ is a word on $X(n)$ and $\y{s_1}{t_1}\cdots\y{s_m}{t_m}$ is a word on $Y(n)$ with the condition that if $s_j \subset s_i$, then $i<j$. 
\end{definition}
We can always make a word in weak standard form into one in standard form. 
\begin{lemma}[{\cite[Lemma 3.11]{lodha2020nonamenable}} for $n$=2]\label{weakstandard_standard}
We can rewrite a weak standard form into a standard form with the same length by just switching the letters $($i.e., relation $(3)$$)$ finitely many times. 
\begin{proof}
Let $f\y{s_1}{t_1}\cdots\y{s_m}{t_m}$ be a word in weak standard form. 
We show this by induction on $m$. 
It is obvious if $m=0, 1$. 
For $m \geq 2$, by the induction hypothesis, we get a word $f\y{s_1^\prime}{t_1^\prime}\cdots\y{s_{m-1}^\prime}{t_{m-1}^\prime}\y{s_m}{t_m}$ where $f\y{s_1^\prime}{t_1^\prime}\cdots\y{s_{m-1}^\prime}{t_{m-1}^\prime}$ is in standard form. 
Then one of the following holds: 
\begin{enumerate}[label=(\alph*)]
\item $s_{m-1}^\prime \supset s_m$; 
\item $s_{m-1}^\prime \perp s_m$ and $s_{m-1}^\prime(i)<s_m(i)$, where $i$ is the smallest number such that $s_{m-1}^\prime (i) \neq s_m(i)$; 
\item $s_{m-1}^\prime \perp s_m$ and $s_{m-1}^\prime(i)>s_m(i)$, where $i$ is the smallest number such that $s_{m-1}^\prime(i)\neq s_m(i)$. 
\end{enumerate}
If (a) or (b), $f\y{s_1^\prime}{t_1^\prime}\cdots\y{s_{m-1}^\prime}{t_{m-1}^\prime}\y{s_m}{t_m}$ is also in standard form. 
If (c), by applying the relation (3) to $\y{s_{m-1}^\prime}{t_{m-1}^\prime}\y{s_m}{t_m}$ and using the induction hypothesis for the first $m-1$ characters again, we get a word in standard form. 
Indeed, let $f\y{s_1^{\prime\prime}}{t_1^{\prime\prime}}\cdots\y{s_{m-1}^{\prime\prime}}{t_{m-1}^{\prime\prime}}\y{s_{m-1}^{\prime}}{t_{m-1}^{\prime}}$ be the obtained form. 
By the transitivity of the order, we only verify that $s_{m-1}^{\prime \prime}<s_{m-1}^{\prime}$ holds. 
There are two cases, $s_{m-1}^{\prime \prime}=s_{j}^{\prime}$ for some $j$ and $s_{m-1}^{\prime \prime}=s_m$.
And in either case, the claim clearly holds. 
\end{proof}
\end{lemma}
For an infinite $n$-ary word $w$ in $\NN$ and a word $f\y{s_1}{t_1}\cdots\y{s_m}{t_m}$ in weak standard form, we define their \textit{calculation} as follows: 
First, we apply $f$ to $w$. Then apply $\y{s_1}{t_1}, \dots, \y{s_m}{t_m}$ to $f(w)$ in this order, where the latter term ``apply'' means to rewrite each $y_{s_i}$ by using its definition in equation \eqref{n-adicy_definition}, and no rewriting be done by the definition of the map $y$. 
\begin{example}
Let $n=4$. For $w=3002\cdots$ and $y_{300}^{-1}y_{30}y_1$, we calculate as follows: 
\begin{align*}
3002\cdots \xrightarrow{y_{300}^{-1}}300y^{-1}(2\cdots) \xrightarrow{y_{30}} 30(y(0(y^{-1}(2\cdots)))) \xrightarrow{y_1} 30(y(0(y^{-1}(2\cdots)))). 
\end{align*}
Therefore the calculation of $w=3002\cdots$ and $y_{300}^{-1}y_{30}y_1$ is $30(y(0(y^{-1}(2\cdots))))$. 
\end{example}
We write such element as $30y0y^{-1}2\cdots$ and sometimes regard it as a word on $\N \cup \{y, y^{-1}\}$. 
We also define calculations for finite words in $\Nn$ and weak standard forms in the same way, although finite words are not in $\NN$. 
We note that, unlike infinite words, not all calculations of finite words can be defined. 
For example, in $n=4$, the calculation of $w=3002$ and $y_{300}^{-1}y_{30}y_1$ is $30y0y^{-1}2$. 
However, that of $w=3$ and $y_{300}^{-1}y_{30}y_1$ is not defined. 

\begin{figure}[tbp]
	\centering
	\includegraphics[width=120mm]{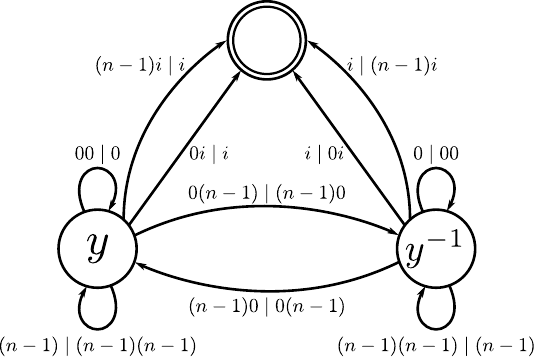}
	\caption{The transducer of $y$, where $i=1, \dots, n-2$. }
	\label{transducer_y}
\end{figure}

We call a \textit{substitution} the operation of rewriting a calculation once using the definition of $y$. 
We also call the given element $w$ and the element obtained by repeating substitutions \textit{input string} and \textit{output string}, respectively. 
\begin{example}
For $w=0^4(n-1)0^4(n-1)\cdots$ and $y^2$, we have the following substitutions: 
\begin{align*}
y^20^4(n-1)0^4(n-1)\cdots \to y 0 y 0^2(n-1)0^4(n-1)\cdots \to y0^2 y(n-1)0^4(n-1)\cdots. 
\end{align*}
Their output string is $0(n-1)^40(n-1)^4 \cdots$. 
\end{example}
As described in \cite[Section 3.6]{lodha2020nonamenable}, the map $y$ can be expressed as a finite state transducer. 
For the sake of simplicity, we assume $n\geq3$. 
Our transducer is illustrated in Figure \ref{transducer_y}. 
In this setting, each edge is labeled by a form $\sigma \mid \tau$, where $\sigma$ and $\tau$ are in $\Nn$. 
The element $\sigma$ represents an input string, and $\tau$ represents an output string. 
The difference from the case of $n=2$ is that there exists an accepting state (double circle mark). 
This difference corresponds to the fact that $y$ vanishes in $y \sigma$ by substitutions if $\sigma$ is not in $\{0, n-1 \}^{\omega}$ nor $\{0, n-1 \}^{<\omega}$. 

\begin{figure}[tbp]
	\centering
	\includegraphics[width=110mm]{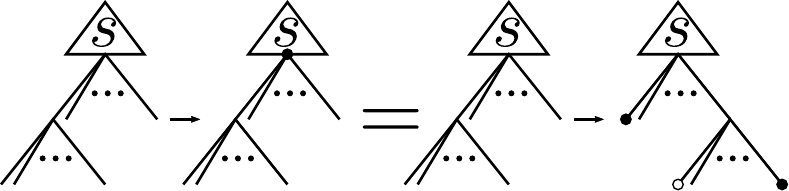}
	\caption{The black and white circles represent that we apply $y$ and $y^{-1}$ to the words corresponding to the edges below them, respectively. }
	\label{image_expansion_move}
\end{figure}

\begin{definition}
A calculation is a \textit{potential cancellation} if this is of the form 
\begin{align*}
y^{t_1}\sigma y^{t_2}, \sigma \in \Nn, t_i \in \{1, -1 \}
\end{align*}
such that we get the word $\sigma^\prime y^{-t_2}$ by substituting $y^{t_1}\sigma$ as many as possible. 
We say a calculation contains a potential cancellation if there exists a subword such that the subword is potential cancellation.  
\end{definition}
\begin{example}
The subwords $y0(n-1)y$ and $y00y^{-1}$ are potential cancellations. 
The subword $y0y$ is not a potential cancellation since we can not apply any substitutions. 
The subword $y^{-1}(n-1)00(n-1)y^{-1}$ also is not a potential cancellation since we have $y^{-1}(n-1)00(n-1)y^{-1} \to 0(n-1)y0(n-1)y^{-1} \to 0(n-1)(n-1)0y^{-1}y^{-1}$. 
\end{example}
\begin{remark} \label{remark_PCA}
If a calculation is a potential cancellation, then $\sigma$ is in $\{0, n-1\}^{<\omega}$. 
Indeed, if not, $y$ vanishes. 
\end{remark}
The following holds. 
\begin{lemma}[{\cite[Lemma 5.9]{lodha2016nonamenable}} for $n=2$] \label{substitution_NPCA}
Let $\Lambda$ be a calculation that contains no potential cancellations. 
Then the word $\Lambda^{\prime}$ obtained by a substitution at any $y^{\pm1}$ again contains no potential cancellations. 
\begin{proof}
By Remark \ref{remark_PCA}, it is sufficient to consider the case $\sigma \in \{0, n-1\}^{<\omega}$. 
Then we can identify $\{0, n-1 \}$ with $\{0, 1\}$. 
By \cite[Lemma 5.9]{lodha2016nonamenable}, we have the desired result. 
\end{proof}
\end{lemma}
We generalize the notion of potential cancellation to the case of weak standard forms. 
\begin{definition}
A weak standard form $f\y{s_1}{t_1}\cdots\y{s_m}{t_m}$ has a \textit{potential cancellation} if there exists $w$ in $\NN$ such that the calculation of $w$ by $f\y{s_1}{t_1}\cdots\y{s_m}{t_m}$ contains a potential cancellation. 
\end{definition}
To construct unique words, we define the following moves obtained from the relations in $R(n)$. 
\begin{definition} \label{definition_five_moves}
We assume that every map in the following is defined. 
\begin{description}
\item[Rearranging move] 
$\y{t}{i}\x{j}{s}^{\delta} \to \x{j}{s}^{\delta}\y{\x{j}{s}^{\delta}(t)}{i}$ ($\delta \in \{+1, -1\}$); 
\item[Expansion move]
$y_s \to \x{0}{s}y_{s0}\y{s(n-1)0}{-1}y_{s(n-1)(n-1)}$ and $y_s^{-1} \to \x{0}{s}^{-1}\y{s00}{-1}y_{s0(n-1)}\y{s(n-1)}{-1}$; 
\item[Commuting move] $y_uy_v\leftrightarrow y_v y_u$ (if $u \perp v$); 
\item[Cancellation move] $\y{s}{\delta}\y{s}{-\delta} \to \epsilon$ ($\delta \in \{+1, -1\}$); 
\item[ER moves] 
\begin{enumerate}
\item $f(\y{s_1}{t_1}\cdots\y{s_k}{t_k})y_u \to f\x{0}{u}(\y{\x{0}{u}(s_1)}{t_1}\cdots \y{\x{0}{u}(s_k)}{t_k})(y_{u0}\y{u(n-1)0}{-1}y_{u(n-1)(n-1)})$
\item $f(\y{s_1}{t_1}\cdots\y{s_k}{t_k})\y{u}{-1} \to f\x{0}{u}^{-1}(\y{\x{0}{u}^{-1}(s_1)}{t_1}\cdots \y{\x{0}{u}^{-1}(s_k)}{t_k})(\y{u00}{-1}y_{u0(n-1)}\y{u(n-1)}{-1})$. 
\end{enumerate}
\end{description}
ER moves are a combination of an expansion move and rearranging moves. 
\end{definition}
By the definition of potential cancellations, we note that if we apply an ER move to a standard form that contains a potential cancellation, then either the resulting word also contains a potential cancellation, or we can apply a cancellation move to the resulting word. 
Like the Lodha--Moore group $G_0$, each element of $G_0(n)$ can be represented by tree diagrams. 
Then the expansion move for $y_s$ is the replacement from the left diagram to the right diagram in Figure \ref{image_expansion_move}. 
See \cite[Section 4]{lodha2016nonamenable} for details of the case $n=2$. 

We introduce the following notion, which makes it easy to check whether the moves are defined. 
\begin{definition}
For $s$ in $\Nn$, we write $\| s\|$ for the length of $s$. 
We define the \textit{depth} of a weak standard form $f\y{s_1}{t_1}\cdots\y{s_m}{t_m}$ to be the integer $\min_{1 \leq i \leq m} \|s_i \|$. 
If a weak standard form is on $X(n)$, we define its depth to be $\infty$. 
\end{definition}
We can make arbitrary words into standard forms. 
The proof is given by a natural generalization of \cite[Lemma 5.4]{lodha2016nonamenable}, where $1$ is replaced by $n-1$. 
\begin{proposition}[{\cite[Lemma 5.4]{lodha2016nonamenable}} for $n=2$] \label{Prop_standard_forms_from_words}
Let $l$ be a natural number and $w$ be a word on $Z(n)$. 
Then we can rewrite $w$ into a standard form whose depth is at least $l$ by moves. 
\end{proposition}
Next, we list three lemmas about ER moves. 
These are useful for getting words in standard form without potential cancellations. 
\begin{lemma}[{\cite[Lemma 3.21]{lodha2020nonamenable}} for $n$=2]\label{ER_move_preserve_wsf}
Let $\Lambda = f(\y{s_1}{t_1}\cdots\y{s_k}{t_k})y_u^{\pm1}(\y{p_1}{q_1}\cdots\y{p_m}{q_m})$ be in weak standard form. 
Then the word obtained from $\Lambda$ by an ER move on $y_u^{\pm1}$ is also in weak standard form. 
\begin{sproof}
We show this only for the case of $y_u$. 
Let 
\begin{align*}
f\x{0}{u}(\y{\x{0}{u}(s_1)}{t_1}\cdots \y{\x{0}{u}(s_k)}{t_k})(y_{u0}\y{u(n-1)0}{-1}y_{u(n-1)(n-1)})(\y{p_1}{q_1}\cdots\y{p_m}{q_m})
\end{align*}
be the word obtained from $\Lambda$ by an ER move. 
We show this by considering whether each $s_i$ is independent with $u$ or not. 
If $s_i \perp u$ holds, then $\x{0}{u}(s_i)=s_i$, and there is nothing to do. 
If $s_i \supset u$ holds, since $\x{0}{u}(s_i)$ is defined, $\x{0}{u}(s_i) \supset u$ holds, and it can never be a proper prefix of $u0$, $u(n-1)0$, or $u(n-1)(n-1)$. Since $p_j$ satisfies either $p_j\perp u$ or $p_j \subseteq u$, the claim holds in both cases. 
\end{sproof}
\end{lemma}
\begin{lemma}[{\cite[Lemma 3.22]{lodha2020nonamenable}} for $n$=2]\label{lem_ERmove_preserve_NPCA}
Let $f \lambda_1$ be in weak standard form with no potential cancellations and 
$g\lambda_2$ be a word obtained from $f \lambda_1$ by an ER move. 
Then $g\lambda_2$ is in weak standard form with no potential cancellations. 
\begin{proof}
Let $\lambda_1=\y{s_1}{t_1}\cdots\y{s_m}{t_m}$ with every $t_i$ in $\{\pm1\}$. 
We apply ER move on $\y{u_i}{t_i}$. 
Since $f$ does not affect the existence of potential cancellations, we assume that $f$ is the empty word. 
We only consider the case $t_i=1$ since the cases $t_i=1$ and $t_i=-1$ are shown similarly. 
Let $s_j^{\prime}\coloneqq \x{0}{s_i}(s_j)$ $(1 \leq j \leq i-1)$. 
By the definition of ER move, we have
\begin{align*}
g\lambda_2=\x{0}{s_i}(\y{s_1^\prime}{t_1}\cdots\y{s_{i-1}^\prime}{t_{i-1}})(y_{s_i0}\y{s_i(n-1)0}{-1}y_{s_i(n-1)(n-1)})(\y{s_{i+1}}{t_{i+1}}\cdots\y{s_m}{t_m}). 
\end{align*}

Since it is clear by Lemma \ref{ER_move_preserve_wsf} that $g\lambda_2$ is in weak standard form, suppose that $g\lambda_2$ has a potential cancellation. 
Let $\tau \in \NN$ be an element such that the calculation $\Lambda$ of $\tau$ by $(\y{s_1^\prime}{t_1}\cdots\y{s_{i-1}^\prime}{t_{i-1}})(y_{s_i0}\y{s_i(n-1)0}{-1}y_{s_i(n-1)(n-1)})(\y{s_{i+1}}{t_{i+1}}\cdots\y{s_m}{t_m})$ contains a potential cancellation. 
Then, $s_i$ is a prefix of $\tau$. 
Indeed, we have the following: 
\begin{enumerate}
\item if $s_j \perp s_i$, then $\y{s_j^\prime}{t_j}=\y{s_j}{t_j}$; 
\item if $s_j \perp s_i$, then $\y{s_j}{t_j}(\tau^\prime)\perp s_i$ where $\tau^\prime$ is in $\NN$ with $\tau^\prime \perp s_i$; 
\item if $s_j \supset s_i$, since $\x{0}{s_i}(s_j) \supset s_i$, we have $\y{s_j^\prime}{t_j}(\tau^\prime)=\y{\x{0}{s_i}(s_j)}{t_j}(\tau^\prime)=\tau^\prime$ and $y_{s_j}(\tau^\prime)=\tau^\prime$ where $\tau^\prime$ is in $\NN$ with $\tau^\prime \perp s_i$; 
\item $y_{s_i}$, $y_{s_i0}$, $y_{s_i(n-1)0}^{-1}$, and $y_{s_i(n-1)(n-1)}$ all fix $\tau^\prime$ where $\tau^\prime$ is in $\NN$ with $\tau^\prime \perp s_i$. 
\end{enumerate}
Since ER move is defined, each $y_j$ satisfies either $s_j\perp s_i$ or $s_j \supset s_i$. 
If $s_i$ is not a prefix of $\tau$, then $\tau \perp s_i$ holds. 
Then the calculations of $\tau$ with $f\lambda_1$ and with $g\lambda_2$ are the same. 
This contradicts the assumption that $f\lambda_1$ does not have a potential cancellation. Hence $s_i$ is a prefix of $\tau$. 

Let $\Lambda^\prime$ be the calculation of $\x{0}{s_i}^{-1}(\tau)$ by $\y{s_1}{t_1}\cdots\y{s_m}{t_m}$. 
By the assumption, this calculation does not contain a potential cancellation. 
On the other hand, $\Lambda$ is obtained from $\Lambda^\prime$ by a single substitution. 
However, this also contradicts Lemma \ref{substitution_NPCA}. 
\end{proof}
\end{lemma}
\begin{lemma} \label{lem_depth_ER_moves}
Let $l$ be an arbitrary natural number. 
Then, by applying ER moves to a weak standard form $f\y{s_1}{t_1}\cdots\y{s_m}{t_m}$ finitely many times, we obtain a weak standard form whose depth is at least $l$. 
\begin{sproof}
If $\|s_i\|< l$ holds, we apply ER moves on $\y{s_i}{t_i}$ (if $t_i \neq \pm1$, we apply on the last one). 
We note that the move may not be defined. 
In that case, before we do that, we apply ER moves to $\y{s_j}{t_j}$ $(j<i)$ which is the cause that the ER move on $\y{s_i}{t_i}$ is not defined. 
If it also is not defined, repeat the process. 
Since we can always apply ER moves on $f\y{s_1}{t_1}$, this process terminates. 
We do this process repeatedly until the depth of the obtained word is at least $l$. 
By Lemma \ref{ER_move_preserve_wsf}, it is also in weak standard form. 
\end{sproof}
\end{lemma}

Since subwords play an essential role in the notion of potential cancellations, we introduce the following definitions for simplicity. 
\begin{definition}
Let $f\y{s_1}{t_1}\cdots\y{s_m}{t_m}$ be in weak standard form. 
We say the pair $(\y{s_j}{t_j}, \y{s_i}{t_i})$ is \textit{adjacent} if the following two conditions hold: 
\begin{enumerate}
\item $s_i \subset s_j$, 
\item if $u$ in $\Nn$ satisfies $s_i \subset u \subset s_j$, then $u \notin \{ s_1, \dots, s_m\}$. 
\end{enumerate}
\end{definition}
We define an adjacent pair $(\y{s_j}{t_j}, \y{s_i}{t_i})$ is a \textit{potential cancellation} if $y^{t_i}\sigma y^{t_j}$ is a potential cancellation, where $\sigma$ is the word satisfying $s_i\sigma=s_j$. 
It is clear from the definition that a weak standard form contains a potential cancellation if and only if there exists an adjacent pair such that it is a potential cancellation. 

By the definition, for example, if $(\y{j}{t}, y_s)$ is a potential cancellation, then either $\y{j}{t}=\y{s00}{-1}$, $y_{0(n-1)}$, or $\y{s(n-1)}{-1}$ holds, or $(\y{\x{0}{s}(j)}{t}, y_{s0})$, $(\y{\x{0}{s}(j)}{t}, \y{s(n-1)0}{-1})$, or $(\y{\x{0}{s}(j)}{t}, y_{s(n-1)(n-1)})$ is also a potential cancellation. 
As mentioned in Remark \ref{remark_PCA}, we note that if an adjacent pair is a potential cancellation, then $\sigma$ is in $\{0, n-1\}^{<\omega}$.

We introduce the moves which are ``inverse moves'' of ER moves. 
First, we define the conditions where the moves are defined. 
\begin{definition}\label{def_potential_contraction}
We say that a weak standard form contains a \textit{potential contraction} if it satisfies either of the following: 
\begin{enumerate}
\item there exists a subword $y_{s0}\y{s(n-1)0}{-1}y_{s(n-1)(n-1)}$, but there does not exist $\y{s(n-1)}{\pm1}$; 
\item there exists a subword $\y{s00}{-1}y_{s0(n-1)}\y{s(n-1)}{-1}$, but there does not exist $\y{s0}{\pm1}$. 
\end{enumerate}
We also call the same when the condition above is satisfied by commuting moves. 
\end{definition}
\begin{example} \label{example_potential_contraction}
Let $n=4$. A weak standard form $y_{300}y_{3030}^{-1}y_{3031}y_{3033}y_{30}$ contains a potential contraction. 
\end{example}
We now define moves. 
\begin{definition}\label{def_contraction_move}
Let $f\lambda$ be in standard form. 
We assume that this form contains a potential contraction in the sense of (1). 
We apply commuting move to all $\y{u}{v}$ ($u \in \{v \mid s0<v \leq s(n-1)(n-1)\}$) except the subword in the standard form to the left of $y_{s0}$ while preserving the order, and we replace $y_{s0}\y{s(n-1)0}{-1}y_{s(n-1)(n-1)}$ with $\x{0}{s}^{-1}y_s$. 
Then, move $\x{0}{s}^{-1}$ to just next to $f$ by rearranging moves. 
Finally, we apply cancellation moves if necessary. 
Note that the rearranging moves are defined, and the word obtained by the above sequence of moves is also in standard form (Lemma \ref{lem_contraction_move}). 
We call this sequence of moves a \textit{contraction move}. 
We also define a contraction move for (2) in a similar way. 
\end{definition}
The moves and the claim in this definition are justified by the following lemma. 
We describe only case (1), but the similarity holds for case (2). 
\begin{lemma}\label{lem_contraction_move}
In Definition $\ref{def_contraction_move}$, the rearranging moves are defined, and the resulting words are again in standard form. 
\begin{proof}
By the definition of standard form, for $y_u^v$, we have either $s0\perp u$ or $s0 \subseteq u$. 
In the former case, we have either $s \perp u$ or $s \subset u$. 
In the latter case, we have $s \subset u$. 
By the definition of potential contraction, there does not exist $\y{s(n-1)}{\pm1}$. 
Hence $\x{0}{s}^{-1}(u)$ is defined, and the rearranging moves can be done. 

Suppose that we have finished applying a contraction move. 
Let $s_j^\prime \coloneqq \x{0}{s}^{-1}(s_j)$, and $f\x{0}{s}^{-1}(\y{s_1^\prime}{t_1}\cdots\y{s_k^\prime}{t_k})y_s(\y{p_1}{q_1}\cdots\y{p_m}{q_m})$ be the resulting word. 
If $s_k^\prime=s_k$, then $s_k^\prime\perp s$ holds. Since the original word is in standard form, we have $s_k^\prime <s$. 

If $s_k^\prime\neq s_k$, then $s_k^\prime \supset s$. 
By the transitivity of this order, $\x{0}{s}^{-1}(\y{s_1^\prime}{t_1}\cdots\y{s_k^\prime}{t_k})y_s$ is in standard form. 

We show that $y_s(\y{p_1}{q_1}\cdots\y{p_m}{q_m})$ is in standard form. 
We note that $s(n-1)(n-1)<p_1$ holds. 
If $s(n-1)(n-1) \supset p_1$, either $s=p_1$ or $s \supset p_1$ holds. 
In the former case, by $y_s\y{p_1}{t_1}\to\y{s}{t_1+1}$ or $\to \epsilon$, we get the standard form. 
In the latter case, we have $s<p_1$. 
If $s(n-1)(n-1) \perp p_1$, since $(n-1)$ is the largest number in $\{0, \dots, n-1\}$, we have $s<p_1$. 
By the transitivity, $\x{0}{s}^{-1}(\y{s_1^\prime}{t_1}\cdots\y{s_k^\prime}{t_k})y_s(\y{p_1}{q_1}\cdots\y{p_m}{q_m})$ is in standard form. 
\end{proof}
\begin{example}
In example \ref{example_potential_contraction}, we obtain $\x{0}{30}^{-1}y_{301}\y{30}{2}$ by applying a contraction move. 
\end{example}
\end{lemma}
Contraction moves have the following property: 
\begin{lemma}\label{lem_contraction_move_preserve_NPCA}
Let $f\y{s_1}{t_1}\cdots\y{s_m}{t_m}$ be in standard form which contains a potential contraction and no potential cancellations. 
Then the word obtained by contraction move contains no potential cancellations. 
\begin{proof}
We show only the case of potential contraction condition (1).  
We assume that we have produced a potential cancellation after applying a contraction move. 
If $\y{s}{-1}$ is in $f\y{s_1}{t_1}\cdots\y{s_m}{t_m}$, then $(y_{s0}, \y{s}{-1})$ is a potential cancellation. 
Therefore we can assume that $\y{s}{k}$ $(k>0)$ is in the obtained word. 
If there exists an adjacent pair that is potential cancellation, then it is in the form of one of $(\y{s}{k}, \y{s^{\prime}}{k^{\prime}})$, $(\y{s^{\prime \prime}}{k^{\prime \prime}}, \y{s}{k})$, or $(\y{u_1}{v_1}, \y{u_2}{v_2})$, where $u_1, u_2 \neq s$. 
In all cases, it contradicts that $f\y{s_1}{t_1}\cdots\y{s_m}{t_m}$ contains no potential cancellations. 
\end{proof}
\end{lemma}
Our goal in the rest of this section is to give the unique word for each element in $G_0(n)$, which satisfies the following: 
\begin{definition} \label{normal_form}
Let $f\y{s_1}{t_1}\cdots\y{s_m}{t_m}$ be in standard form with no potential cancellations, no potential contractions, and $f$ is in the normal form in the sense of $F(n)$. 
Then we say that $f\y{s_1}{t_1}\cdots\y{s_m}{t_m}$ is in the \textit{normal form}. 
\end{definition}
See Section \ref{subsection_F(n)_definition} for the definition of the normal form of elements in $F(n)$. 

We obtain a form in normal form from an arbitrary word on $Z(n)$ through the following four steps: 
\begin{description}
\item[Step 1] Convert an arbitrary word into a standard form (Proposition \ref{Prop_standard_forms_from_words}); 
\item[Step 2] Convert the word into a standard form that contains no potential cancellations; 
\item[Step 3]Convert the word into a standard form $f\y{s_1}{t_1}\cdots\y{s_m}{t_m}$ which contains no potential cancellations and no potential contractions; 
\item[Step 4]Convert $f$ into $g$, where $g$ is the unique normal form in $F(n)$ (Section \ref{subsection_F(n)_definition}). 
\end{description}
The remaining steps are only 2 and 3. 
\begin{lemma}[{\cite[Lemma 4.5]{lodha2020nonamenable}} for $n$=2]\label{lem_for_step2}
Let $(\y{s_1}{t_1}\cdots\y{s_m}{t_m})\y{s}{t}$ be in standard form such that the following hold: 
\begin{enumerate}[font=\normalfont]
\item $t_1, \dots, t_m, t \in \{1, -1\}$; 
\item any two $s_i, s_j$ $(i \neq j)$ are independent; 
\item for any $\y{s_i}{t_i}$, $(\y{s_i}{t_i}, \y{s}{t})$ is an adjacent pair with being potential cancellation. 
\end{enumerate}
Then by applying ER moves and cancellation moves, we obtain a standard form $f\y{u_1}{v_1}\cdots\y{u_k}{v_k}$ such that any two $u_i, u_j$ are independent and $v_1, \dots, v_k \in \{1, -1\}$. 
\begin{proof}
We consider only the case $t=1$. 
We claim that $\x{0}{s}(s_i)$ is defined for every $i$. 
Indeed, for $s_i$ such that $s \subset s_i$, the only case where $\x{0}{s}(s_i)$ is not defined is the case $s_i=s0$. 
Since the adjacent pair $(\y{s0}{t_i}, y_s)$ is not potential cancellation whether $t_i$ is $1$ or $-1$, this contradicts the assumption $(3)$. 

Let $s_i^\prime=\x{0}{s}(s_i)$. 
By applying an ER move to $y_s$, we have
\begin{align*}
(\y{s_1}{t_1}\cdots\y{s_m}{t_m})y_s=\x{0}{s}(\y{s_1^\prime}{t_1}\cdots\y{s_m^\prime}{t_m})y_{s0}\y{s(n-1)0}{-1}y_{s(n-1)(n-1)}. 
\end{align*}
We note that $s0$, $s(n-1)0$, and $s(n-1)(n-1)$ are independent of each other. 
For every $\y{s_i^\prime}{t_i}$, one of $y_{s0}$, $\y{s(n-1)0}{-1}$, or $y_{s(n-1)(n-1)}$ is its inverse element or forms an adjacent pair.  
Let $\y{\sigma_i}{\tau_i}$ be the corresponding one. 
If it is the inverse, we apply a cancellation move. 
If it forms an adjacent pair, the distance in an $n$-ary tree between $s_i^\prime$ and $\sigma$ is smaller than the distance between $s_i$ and $s$. 
Thus, by iterating this process, we obtain the desired result. 
\end{proof}
\end{lemma}
The following lemma completes step 2. 
\begin{lemma}[{\cite[Lemma 4.6]{lodha2020nonamenable}} for $n$=2]\label{lem_step2}
By applying moves, any weak standard form can be rewritten into a standard form that contains no potential cancellations. 
\begin{proof}
Let $f\y{s_1}{t_1}\cdots\y{s_m}{t_m}$ $(t_i \in \{-1, 1 \})$ be a weak standard form. 
We show by induction on $m$. 
If $m\leq1$, since there exists no adjacent pair, it is clear. 

We assume that $m>1$ holds. 
By applying the inductive hypothesis, Lemmas \ref{lem_ERmove_preserve_NPCA}, and \ref{lem_depth_ER_moves} to $f\y{s_1}{t_1}\cdots\y{s_{m-1}}{t_{m-1}}$, we obtain a weak standard form $g\y{p_1}{q_1}\cdots \y{p_k}{q_k}$ with at least depth $\|s_m\|+1$ and without potential cancellations. 
Since $p_i \subset s_m$ does not hold by its depth, $g\y{p_1}{q_1}\cdots \y{p_k}{q_k}\y{s_m}{t_m}$ is in weak standard form. 

If there exists an adjacent pair that is a potential cancellation, it is only of the form $(\y{p_i}{q_i}, \y{s_m}{t_m})$. 
We record all such $\y{p_i}{q_i}$. 
By applying commuting moves, we obtain a weak standard form $h(\y{j_1}{k_1}\cdots\y{j_l}{k_l})(\y{u_1}{v_1}\cdots\y{u_o}{v_o})\y{s_m}{t_m}$ which satisfies the following: 
\begin{enumerate}
\item $v_1, \dots, v_o \in \{ -1, 1\}$; 
\item For $i=1, \dots, o$, $(\y{u_i}{v_i}, \y{s_m}{t_m})$ is an adjacent pair which is a potential cancellation; 
\item All other adjacent pairs are not potential cancellations. 
\end{enumerate}
Indeed, since $s_m$ is the shortest word, each adjacent element is the ``second shortest.'' 
Hence we can apply commuting moves to get a word. 

By Lemmas \ref{lem_ERmove_preserve_NPCA} and \ref{lem_depth_ER_moves}, note that we can increase the depth of the word $h(\y{j_1}{k_1}\cdots\y{j_l}{k_l})$ by ER moves while preserving it in weak standard form without potential cancellations. 
Hence we can assume that the depth of $h(\y{j_1}{k_1}\cdots\y{j_l}{k_l})$ is sufficiently large. 
Then we  apply ER moves to $h(\y{j_1}{k_1}\cdots\y{j_l}{k_l})(\y{u_1}{v_1}\cdots\y{u_o}{v_o})\y{s_m}{t_m}$ as in Lemma \ref{lem_for_step2}. 
By the same argument as in Lemma \ref{lem_ERmove_preserve_NPCA}, no new potential cancellations are produced in this process. 
Finally, by Lemma \ref{weakstandard_standard}, we have the desired result. 
\end{proof}
\end{lemma}
The following lemma completes step 3. 
\begin{lemma}[{\cite[Lemma 4.7]{lodha2020nonamenable}} for $n$=2]\label{lem_step3}
Any standard form which contains no potential cancellations can be rewritten into a standard form that contains no potential cancellations and no potential contractions. 
\begin{proof}
We apply contraction moves repeatedly. 
By Lemmas \ref{lem_contraction_move} and \ref{lem_contraction_move_preserve_NPCA}, the applied word is again in standard form and contains no potential cancellations. 
Since this move makes the word on $Y(n)$ in the standard form strictly shorter, the process finishes. 
\end{proof}
\end{lemma}
The only thing that remains to be proved is the uniqueness of the normal form. 
We will show this by contradiction. 
In order to do this, we define an ``invariant'' of the forms. 
\begin{definition}
A calculation that contains no potential cancellations has \textit{exponent $m$} if $m$ is the number of $y^{\pm1}$ satisfying the condition that the letter appearing after it is only $0$, $(n-1)$, or $y^{\pm1}$ and no other letters appear. 
\end{definition}
\begin{example}
Let $n=4$. 
The string $y0y02y^{-1}3y00\cdots$ has exponent $2$ since there exists 2. 
The string $y100\cdots$ has exponent $0$ since there exists $1$. 
\end{example}
For the proof of uniqueness, we prepare two lemmas. 
\begin{lemma}[{\cite[Lemma 5.10]{lodha2016nonamenable}} for $n=2$]\label{lem_exponent_1}
Let $\Lambda$ be a finite word on $\N \cup \{y, y^{-1}\}$ that contains no potential cancellations. 
Let $m$ be the exponent of $\Lambda$. 
Then there exists a finite word $u$ on $\{0, n-1\}$ and $v$ in $\Nn$ such that $\Lambda u$ can be rewritten into $vy^m$ by substitutions. 
\begin{proof}
By substitutions, we rewrite $\Lambda$ into $v^\prime\Lambda^\prime$, where $v^\prime$ is in $\Nn$ and $\Lambda^\prime$ is a word on $\{0, (n-1), y, y^{-1}\}$. 
By the definition of exponent, $\Lambda^\prime$ also has exponent $m$. 
Also, by Lemma \ref{substitution_NPCA}, this contains no potential cancellations. 
Hence, by identifying $\{0, n-1\}$ with $\{0, 1\}$, the claim comes down to the case $n=2$. 
\end{proof}
\end{lemma}
For $u, s$ in $\Nn$, we say that \textit{$u$ dominates $s$} if the condition $u \perp s$ or $u \supset s$ is satisfied. 
Then the following holds. 
The proof is given by replacing $1$ with $n-1$ of \cite[Lemma 4.8]{lodha2020nonamenable}. 
\begin{lemma}[{\cite[Lemma 4.8]{lodha2020nonamenable}} for $n=2$]\label{lem_exponent_2}
Let $f\y{s_1}{t_1}\cdots\y{s_l}{t_l}$ and $g\y{p_1}{q_1}\cdots\y{p_m}{q_m}$ be in standard form which represent the same element in $G_0(n)$. 
Let $u$ be in $\Nn$ such that the following hold: 
\begin{enumerate}[font=\normalfont]
\item $f(u)\eqqcolon u_1$ and $g(u)\eqqcolon u_2$ are defined; 
\item $u_1$ dominates $s_1, \dots, s_l$; 
\item $u_2$ dominates $p_1, \dots, p_m$. 
\end{enumerate}
Let $\Theta$ be the calculation of $u$ and $f\y{s_1}{t_1}\cdots\y{s_l}{t_l}$, and let $\Lambda$ be the calculation of $u$ and $g\y{p_1}{q_1}\cdots\y{p_m}{q_m}$. 
Assume that $\Theta$ and $\Lambda$ contain no potential cancellations. 
Then their exponents are the same. 
\end{lemma}
\begin{corollary} \label{G_0(n)_presentation}
The group obtained from the presentation $\langle Z(n) \mid R(n)\rangle$ is the $n$-adic Lodha--Moore group $G_0(n)$. 
\begin{proof}
We show that any standard form which contains no potential cancellations and represents an element of $F(n)$ is always a word on $X(n)$. 
It is clear that if $\y{s_1}{t_1}\cdots\y{s_m}{t_m}$ ($m \geq 1$) is in standard form with no potential cancellations, then the exponent of $s_100\cdots \in \NN$ and $\y{s_1}{t_1}\cdots\y{s_m}{t_m}$ is greater than zero. 
Since every exponent of any word in $F(n)$ and any elements in $\NN$ is zero, by Lemma \ref{lem_exponent_2}, we have a desired result. 
%
\end{proof}
\end{corollary}
\begin{remark}
From the above argument, the uniqueness of the normal form of the elements of $F(n)$ in $G_0(n)$ follows. 
\end{remark}\label{remark_uniqueF(n)}
The following completes the proof of the uniqueness. 
Note that our proof is a slightly modified version of \cite{lodha2020nonamenable}. 
\begin{theorem}[{\cite[Theorem 4.4]{lodha2020nonamenable}} for $n=2$]\label{thm_normal_form_uniqueness}
For each element in $G_0(n)$, its normal form is unique. 
\begin{proof}
By Remark \ref{remark_uniqueF(n)}, we can assume that every normal form in the following argument is not in $F(n)$. 

We show this by contradiction. 
Let $f\y{s_1}{t_1}\cdots\y{s_l}{t_l}$ and $g\y{p_1}{q_1}\cdots\y{p_m}{q_m}$ be different normal forms representing the same element. 
We can assume that $s_l \leq p_m$ without loss of generality. 
One of the following three holds: 
\begin{enumerate}
\item $s_l=p_m$ and $t_l=q_m$; 
\item $s_l=p_m$ and $t_l \neq q_m$; 
\item $s_l<p_m$. 
\end{enumerate}
First, we show that it is sufficient to consider only case $(3)$. 

In case (1), since $\y{s_l}{t_l}=\y{p_m}{q_m}$, we start from $f\y{s_1}{t_1}\cdots\y{s_{l-1}}{t_{l-1}}$ and $g\y{p_1}{q_1}\cdots\y{p_{m-1}}{q_{m-1}}$. 
They are two different standard forms representing the same element. 
They contain no potential cancellations but possibly contain a potential contraction. 
Since we only cancel $\y{s_l}{t_l}=\y{p_m}{q_m}$, the only case that contains a potential cancellation at this step is when $\y{s_l}{t_l}$ or $\y{p_m}{q_m}$ or both play the role of $\y{s(n-1)}{\pm1}$ for some $s$ (see Definition \ref{def_potential_contraction} (1)). 
By the assumption of the standard form, the case with $\y{s0}{\pm 1}$ (Definition \ref{def_potential_contraction} (2)) does not occur. 
We further divide case (1) into the following four subcases: 
\begin{description}
\item[\rm (1-1)] Both $f\y{s_1}{t_1}\cdots\y{s_{l-1}}{t_{l-1}}$ and $g\y{p_1}{q_1}\cdots\y{p_{m-1}}{q_{m-1}}$ contain a potential contraction, and $t_{l-1}=q_{m-1}$; 
\item[\rm (1-2)] Both $f\y{s_1}{t_1}\cdots\y{s_{l-1}}{t_{l-1}}$ and $g\y{p_1}{q_1}\cdots\y{p_{m-1}}{q_{m-1}}$ contain a potential contraction, and $t_{l-1}\neq q_{m-1}$; 
\item[\rm (1-3)] $f\y{s_1}{t_1}\cdots\y{s_{l-1}}{t_{l-1}}$ contains a potential contraction, and $g\y{p_1}{q_1}\cdots\y{p_{m-1}}{q_{m-1}}$ contains no potential contractions, or vice versa; 
\item[\rm (1-4)] $f\y{s_1}{t_1}\cdots\y{s_{l-1}}{t_{l-1}}$ and $g\y{p_1}{q_1}\cdots\y{p_{m-1}}{q_{m-1}}$ contain no potential contractions. 
\end{description}

In case (1-1), we consider $f\y{s_1}{t_1}\cdots\y{s_{l-2}}{t_{l-2}}$ and $g\y{p_1}{q_1}\cdots\y{p_{m-2}}{q_{m-2}}$ instead of the original words in order to eliminate the potential contraction part. 
Since the assumption of containing a potential contraction and being in standard form, $s_l=p_m=s(n-1)$, $s_{l-1}=p_{m-1}=s(n-1)(n-1)$ and $t_{l-1}=q_{m-1}>0$ hold. 
Then $f\y{s_1}{t_1}\cdots\y{s_{l-2}}{t_{l-2}}$ and $g\y{p_1}{q_1}\cdots\y{p_{m-2}}{q_{m-2}}$ are different normal forms. 
Indeed, the only case where a potential contraction is generated by canceling $\y{s_{l-1}}{t_{l-1}}=\y{p_{m-1}}{q_{m-1}}$ is when $\y{s_{l-1}}{t_{l-1}}=\y{p_{m-1}}{q_{m-1}}$ plays the role of $y_{s^\prime(n-1)}$ for $s^\prime=s(n-1)$, but this does not occur since $\y{s^\prime0}{k_1}=\y{s(n-1)0}{k_1}$ and $\y{s^\prime0}{k_2}=\y{s(n-1)0}{k_2}$ exist in two forms respectively and $k_1, k_2<0$ holds by the assumption of potential contractions. 
Thus, we get ``shorter'' normal forms that satisfy all the assumptions. 

In case (1-2), we can assume that $0<t_{l-1}<q_{m-1}$ without loss of generality. 
We consider $f\y{s_1}{t_1}\cdots\y{s_{l-2}}{t_{l-2}}$ and $g\y{p_1}{q_1}\cdots\y{p_{m-1}}{q_{m-1}-t_{l-1}}$ instead. 
For the same reason as for case (1-1), the former is in normal form. 
By the assumption of a potential contraction, $\y{s_{l-2}}{t_{l-2}}=\y{s(n-1)i\sigma}{t_{l-2}}$, where $i$ is in $\N$ and $\sigma$ is a word in $\Nn$. 
Note that $s(n-1)0 \leq s(n-1)i \sigma$ holds. 
For the latter, by performing contraction moves as in step 3, the last letter is $\y{u}{k}$ where $u \subseteq s$. 
This is the situation in case (3). 

In case (1-3), we only consider the former situation. 
By the assumptions, $s_l=p_m=s(n-1)$ holds for some $s$. 
By applying contraction moves to $f\y{s_1}{t_1}\cdots\y{s_{l-1}}{t_{l-1}}$ as in step 3, the last letter is $\y{u}{k}$ where $u \subseteq s \subset s_l=p_m$. 
Since $g\y{p_1}{q_1}\cdots\y{p_{m-1}}{q_{m-1}}$ is in normal form, and in particular in standard form, we have $y_{p_{m-1}}\neq y_u$. 
This is the situation in case (3). 

In case (1-4), both are ``shorter'' normal forms that satisfy all the assumptions. 
Therefore, any subcase of case (1) yields shorter words or case (3). 

In case (2), we can assume that neither $0<q_m<t_l$ nor $t_l<q_m<0$ holds without loss of generality. 
Consider $f\y{s_1}{t_1}\cdots\y{s_{l-1}}{t_{l-1}}$ and $g\y{p_1}{q_1}\cdots\y{p_{m}}{q_{m}-t_{l}}$ instead. 
The latter is in normal form, but as in case (1), the former may contain a potential contraction. 
Similarly, we divide case (2) into two subcases: 
\begin{description}
\item[\rm (2-1)] $f\y{s_1}{t_1}\cdots\y{s_{l-1}}{t_{l-1}}$ contains a potential contraction; 
\item[\rm (2-2)] $f\y{s_1}{t_1}\cdots\y{s_{l-1}}{t_{l-1}}$ contains no potential contractions. 
\end{description}
In case (2-1), the same argument as in (1-3) can be applied. 
In case (2-2), since $s_{l-1}\neq p_m$, this is the situation in case (3). 

We note that case (1-1) or (1-4) happens only finitely many times due to the uniqueness of the normal form of $F(n)$. 
Thus, we consider case (3). 

We only consider the case $q_m>0$. 
We note that $f\y{s_1}{t_1}\cdots\y{s_l}{t_l}\y{p_m}{-1}$ is in standard form since $s_l < p_m$. 
One of the following holds: 
\begin{enumerate}[label=(\roman*)]
\item there exists an infinite word $\sigma$ on $\{0, n-1\}$ such that for any finite word $\sigma_1 \subset \sigma$, $p_m\sigma_1$ is not in $\{s_1, \dots, s_l \}$; 
\item there exists an adjacent pair of the form $(s_i, p_m)$ where $p_mu=s_i$ for $u$ on $\{0, n-1\}$ which is not a potential cancellation. 
\end{enumerate}
Indeed, if both are false, it contradicts that $f\y{s_1}{t_1}\cdots\y{s_l}{t_l}$ contains no potential contractions. 
Then, in either case, there exists a finite word $w$ on $\{0, n-1\}$ such that the calculation $\Lambda$ of $f\y{s_1}{t_1}\cdots\y{s_l}{t_l}\y{p_m}{-1}$ and $f^{-1}(p_mw)$ contains no potential cancellations. 
By considering $w00\cdots 0$ instead,  if necessary, we can consider the following three calculations: 
\begin{enumerate}[label=(\alph*)]
\item $\Theta$ is the calculation of $g\y{p_1}{q_1}\cdots\y{p_m}{q_m-1}$ and $f^{-1}(p_mw)$; 
\item $\Lambda^\prime$ is the calculation of $f\y{s_1}{t_1}\cdots\y{s_l}{t_l}$ and $f^{-1}(p_mw)$; 
\item $\Theta^\prime$ is the calculation of $g\y{p_1}{q_1}\cdots\y{p_m}{q_m}$ and $f^{-1}(p_mw)$. 
\end{enumerate}

All calculations contain no potential cancellations. 
Indeed, for $\Lambda^\prime$ and $\Theta^\prime$, it is clear from the assumption of normal form. 
For $\Theta$, either $q_m-1$ is zero or strictly greater than zero since $q_m>0$. 
In both cases, this is clear from the definitions and the assumption that $g\y{p_1}{q_1}\cdots\y{p_m}{q_m}$ contains no potential cancellations. 

Let $e(\Theta)$, $e(\Theta^\prime)$, $e(\Lambda)$, and $e(\Lambda^\prime)$ be their exponents, respectively. 
Since we have 
\begin{align*}
f\y{s_1}{t_1}\cdots\y{s_l}{t_l}\y{p_m}{-1}=g\y{p_1}{q_1}\cdots\y{p_m}{q_m-1}
\end{align*}
as elements of $G_0(n)$, by Lemma \ref{lem_exponent_2}, $e(\Lambda)=e(\Theta)$ holds. 
Similarly, we have $e(\Theta^\prime)=e(\Lambda^\prime)$. 
By the construction of $\Lambda$ and $\Lambda^\prime$, we have $e(\Lambda)>e(\Lambda^\prime)$. 
Similarly, we have $e(\Theta)\leq e(\Theta^\prime)$ since $q_m>0$. 
By combining them, we obtain
\begin{align*}
e(\Theta)\leq e(\Theta^\prime)=e(\Lambda^\prime)<e(\Lambda), 
\end{align*}
which is a contradiction. 

If $q_m<0$, we can prove the same for $f\y{s_1}{t_1}\cdots\y{s_l}{t_l}y_{p_m}$ and $g\y{p_1}{q_1}\cdots\y{p_m}{q_m+1}$. 
\end{proof}
\end{theorem}

Since $F(n)$ is finitely presented, by routine modifications of \cite[Theorem 3.3]{lodha2016nonamenable}, we also have the following: 
\begin{theorem} \label{Thm_finitely_presented_G0(n)}
$G_0(n)$ admits a finite presentation. 
\end{theorem}

\section{Several properties of $G_0(n)$}\label{section_G0(n)_properties}
We show that $G_0(n)$ has ``expected'' properties. 
\subsection{Nonamenability of $G_0(n)$}
In this section, we discuss embeddings of the groups $G_0(n)$. 
For nonamenability, we only use the fact that a subgroup of an amenable group is also amenable \cite[Theorem 18.29 (1)]{dructu2018geometric}. 
The idea for the following proposition comes from \cite{burillo2001metrics}. 
\begin{theorem}\label{prop_embedding_G_0(p)G_0(q)}
Let $p, q \geq 2$ and assume that there exists $d$ in $\mathbb{N}$ such that $q-1=d(p-1)$ holds. 
Then there exists an embedding $I_{p, q}\colon G_0(p) \to G_0(q)$. 
Moreover, the equality $I_{p,r}=I_{p, q}I_{q, r}$ holds for the maps $I_{p, q}$, $I_{q, r}$ and $I_{p, r}$ defined for $r\geq q \geq p \geq2$ such that any two satisfy the condition. 
\begin{proof}
For the sake of clarity, we label elements of $G_0(p)$ with ``tildes'' and elements of $G_0(q)$ with ``hats.'' 
We first recall the definition of an embedding $F(p) \to F(q)$ given in \cite[Section 3, example (3)]{burillo2001metrics}. 
This homomorphism is defined by
\begin{align*}
&\tilde{x}_0\mapsto\hat{x}_0 & &\tilde{x}_1 \mapsto \hat{x}_d & &\cdots & &\tilde{x}_{p-2} \mapsto \hat{x}_{d(p-2)}& &\tx{0}{(p-1)} \mapsto \hx{0}{(q-1)}, 
\end{align*}
and extended to the (quasi-isometric) embedding \cite[Theorem 6]{burillo2001metrics}. 
By considering $F(p)$ and $F(q)$ as pairs of $p$-ary trees and $q$-ary trees, this embedding is regarded as a ``$p$-caret replacements.'' 
Indeed, for every pair of $p$-ary trees, inserting $d-1$ edges between every pair of adjacent edges of each $p$-caret corresponds to the embedding. 
See Figure \ref{embedding_F3F7}, for example. 
\begin{figure}[tbp]
	\centering
	\includegraphics[width=100mm]{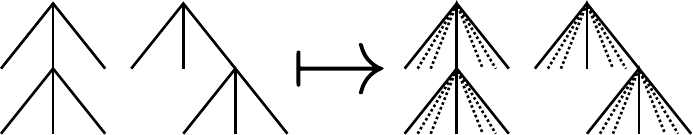}
	\caption{Example of embedding for $p=3$ and $q=7$. }
	\label{embedding_F3F7}
\end{figure}

Define the map $i_{p, q}\colon \Pn=\{0, \dots, p-1\}^{<\omega} \to \Qn=\{0, \dots, q-1\}^{<\omega}$ by setting
\begin{align*}
w_1 w_2\cdots w_k \mapsto (dw_1)(dw_2)\cdots (dw_k), 
\end{align*}
namely, we multiply each letter of a given word by $d$. 
In addition, define the map $Y(p)\to Y(q)$ by setting $\tilde{y}_s \mapsto \hat{y}_{i_{p, q}(s)}$. 
By the definition  of $i_{p, q}$, we have that $\hat{y}_{i_{p, q}(s)}$ is in $Y(q)$. 

By combining the two maps $F(p) \to F(q)$ and $Y(p) \to Y(q)$, we obtain the map $I_{p, q}\colon F(p)\cup Y(p) \to F(q)\cup Y(q)$. 
This map can be extended to the homomorphism $I_{p, q}\colon G_0(p) \to G_0(q)$. 
Indeed, since $i_{p, q}(\tx{i}{s}(t))=\hx{di}{i_{p, q}(s)}(i_{p, q}(t))$ and $I_{p, q}(\tilde{y}_s)=\hat{y}_{i_{p, q}(s)}$ hold, we can verify that the relations of the infinite presentation (Corollary \ref{G_0(n)_presentation}) are preserved under $I_{p, q}$ by direct calculations. 

We claim that $I_{p, q}\colon G_0(p) \to G_0(q)$ is injective. 
We show this by contradiction. 
Assume that $\operatorname{Ker}(I_{p, q})$ is not trivial. 
By the construction, the restriction $I_{p, q}\restr{F(p)}$ coincides with the embedding $F(p) \to F(q)$ mentioned above. 
Hence if there exists $x$ in $\operatorname{Ker}(I_{p, q})$ that is not identity, then $x$ is not in $F(p)$. 
We note that the map $I_{p, q}$ preserves being in normal form. 
In particular, if $x$ is not in $F(p)$, then $I_{p, q}(x)$ is not in $F(q)$ by the uniqueness (Theorem \ref{thm_normal_form_uniqueness}). 
This implies that $I_{p, q}(x)$ is not identity, as desired. 

Let $d_1$, $d_2$ and $d_3$ be natural numbers such that $q-1=d_1(p-1)$, $r-1=d_2(q-1)$ and $r-1=d_3(p-1)$ holds respectively. 
We note that $d_3=d_1d_2$ holds. 
By considering the definitions of the map $i_{p, q}$ and the homomorphism from $F(p)$ to $F(q)$, the equality of $I_{p,r}$ and $I_{p, q}I_{q, r}$ follows immediately. 
\end{proof}
\end{theorem}
\begin{corollary} \label{Cor_G0n_nonamenable}
Let $n\geq2$. Then $G_0(n)$ is nonamenable. 
\begin{proof}
Since we have $n-1=(n-1)(2-1)$ for every $n$, by Theorem \ref{prop_embedding_G_0(p)G_0(q)}, we have an embedding from $G_0(2)=G_0$ into $G_0(n)$. 
This implies that $G_0(n)$ has a subgroup that is isomorphic to $G_0$. 
Since $G_0$ is nonamenable \cite[Theorem 1.1]{lodha2016nonamenable}, $G_0(n)$ is also nonamenable. 
\end{proof}
\end{corollary}
\begin{corollary}
Let $s_i\coloneqq 2^{i-1}+1$. 
Then the sequence $G(s_1), G(s_2), \dots$ forms an inductive system of groups. 
\end{corollary}
\subsection{$G_0(n)$ has no free subgroups}
We assume $n \geq 3$. 
In order to apply the argument of \cite[Theorem 14]{monod2013groups} to $G_0(n)$, we show two lemmas. 

For an element $g$ in $G_0(n)$, we write the set-theoretic support $\supp{g}$ for the set $\{\xi \in \NN \mid g(\xi)\neq\xi\}$. 
Although the space $\NN$ (with the topology) is homeomorphic to the Cantor set, which is totally disconnected, we can consider ``connected components'' of $\supp{g}$ using the total order of $\NN$. 

For $a, b \in \NN$ with $a<b$, set $(a, b) \coloneqq \{\xi \in \NN \mid a<\xi<b\}$. 
\begin{lemma}\label{lem_supp_finite_component}
For $g \in G_0(n)$, there exists a sequence $a_1<a_2 \leq a_3<a_4\leq \cdots < a_{2m}$ in $\NN$ such that we have
\begin{align*}
\supp{g}=(a_1, a_2) \sqcup (a_3, a_4) \sqcup \cdots \sqcup (a_{2m-1}, a_{2m}). 
\end{align*}
We call each $(a_{2i}, a_{2i+1})$ a \textit{connected component} of $\supp{g}$. 
\begin{proof}
Let $f\y{s_1}{t_1}\cdots \y{s_m}{t_m}$ be in normal form of $g$. 
We decompose $\NN$ into disjoint sets $\{n_1\eta \mid \eta \in \NN \}, \dots, \{n_p\eta \mid \eta \in \NN\}$ by using $n_1, \dots, n_p$ in $\Nn$ such that the following holds: 
\begin{enumerate}
\item For each $n_l$, $f(n_l)$ is defined; 
\item For each $s_l$, there exists $n_{l^\prime}$ such that $s_l \leq f(n_{l^\prime})$ holds. 
\end{enumerate}

We show the claim in the lemma by contradiction. 
Assume that there exist $a_1<b_1<\cdots$ such that each $a_i$ is a fixed point, and each $b_i$ is in $\supp{g}$. 
Since $\NN=\{n_1\eta \mid \eta \in \NN \}\cup \dots \cup \{n_p\eta \mid \eta \in \NN\}$ holds, we can assume that each $a_i$ and $b_i$ is in some set $\{n_j\eta \mid \eta \in \NN \}$ without loss of generality. 
Since $a_{i}<b_i<a_{i+1}$ holds, we can write each $a_i$ and $b_i$ as $n_ja_i^\prime$ and $n_jb_i^\prime$, respectively. 

We claim that each $a_i^\prime$ and $b_i^\prime$ can be replaced by words in $\{0, n-1\}^{\omega}$, respectively. 
Indeed, if $a_i^\prime$ is not in $\{0, n-1\}^{\omega}$, then there exist $\hat{a}_i \in \{0, n-1\}^{<\omega}$, $k \in \{1, \dots, n-2\}$, and $a_i^{\prime \prime} \in \NN$ such that $a_i^\prime=\hat{a}_i k a_i^{\prime \prime}$ holds. 
Then, by the definition of $y$, $g(n_j\hat{a}_i k)$ is defined and is equal to $n_j\hat{a}_i k$ since $a_i=n_j \hat{a}_i k a_i^{\prime \prime}$ is a fixed point of $g$. 
Moreover, for any other $k^\prime \in \{1, \dots, n-2\}$, we have $g(n_j\hat{a}_i k^\prime)=n_j\hat{a}_i k^\prime$. 
This implies that we have $g(n_j\hat{a}_i 0 \overline{(n-1)})=n_j\hat{a}_i 0 \overline{(n-1)}$, where $\overline{(n-1)}$ denotes the element $(n-1)(n-1)(n-1)\cdots $ in $\NN$. 
By a similar argument for $b_i^\prime$, we have $g(n_j\hat{b}_i k)=n_j\tilde{b_i}k$, where $\hat{b}_i \neq \tilde{b}_i$. 
Since we have $g(n_j\hat{b}_i k^\prime)\neq n_j\hat{b}_i k^\prime$ for any other $k^\prime$, $n_j\hat{b}_i 0\overline{(n-1)}$ is in $\supp{g}$. 
Even if each $a_i$ or $b_i$ is replaced by the above one, the order $a_1<b_1<\cdots$ is also preserved. 

Let $w_1, \dots, w_t$ be in $\{0, \dots, n-1\}$ such that $f(n_j)=w_1 \cdots w_t$ holds. 
For $n_j \zeta \in \NN$ and $f\y{s_1}{t_1}\cdots \y{s_m}{t_m}$, let $w_1y^{l_1}w_2 y^{l_2}\cdots w_t y^{l_t} \zeta$ be their calculation. 
Then some $l_q$ may be zero. 
Assume that $f(n_j)=w_1 \cdots w_t$ is not in $\{0, n-1\}^{<\omega}$ and let $z$ be the maximal number such that $w_z$ is in $\{1, \dots, n-2\}$. 
By applying finitely many substitutions to all $y^{l_{z^\prime}}$ ($z^\prime<z$), we obtain a word 
$w y^{l_1^\prime} w_2^\prime y^{l_2^\prime} \cdots w_{t^\prime}^\prime y^{l_{t^\prime}^\prime}\zeta$, 
where $w_j^\prime \in \{0,  n-1\}$. 
Since $n_j a_i$ is the output string of $w y^{l_1^\prime} w_2^\prime y^{l_2^\prime} \cdots w_{t^\prime}^\prime y^{l_{t^\prime}^\prime} a_i$, we have $w \leq n_j$. 
Indeed, it is clear that $w \perp n_j$ does not hold, and if $w > n_j$ holds, then it contradicts the facts that $a_i$ is in $\{0, n-1\}^\omega$, and the last number of $w$ is $w_j$. 
Let $n_j=w n_j^\prime$. 
Then we note that $n_j^\prime$ is in $\{0, n-1\}^{<\omega}$ since $y^{l_1^\prime} w_2^\prime y^{l_2^\prime} \cdots w_{t^\prime}^\prime y^{l_{t^\prime}^\prime}a_i$ is in $\{0, n-1\}^\omega$ and its output string equals to $n_j^\prime a_i$. 

We recall that $G_0(n)$ has a subgroup that is isomorphic to $G_0$ (Corollary \ref{Cor_G0n_nonamenable}). 
We construct an element $g^\prime$ which is in the subgroup and satisfies the assumption of $g$. 
Let $g^\prime$ be the element represented by 
\begin{align*}
f^\prime \y{(n-1)0w_2^\prime \cdots w_{t^\prime}^\prime}{l_{t^\prime}^\prime} \y{(n-1)0w_2^\prime \cdots w_{t^\prime-1}^\prime}{l_{t^\prime-1}^\prime} \cdots \y{(n-1)0w_2^\prime}{l_2^\prime} \y{(n-1)0}{l_1^\prime}, 
\end{align*}
where $f^\prime$ is a word on $\{x_0, \x{0}{(n-1)}\}$ satisfying $f^\prime((n-1)0n_j^\prime)=(n-1)0w_2^\prime \cdots w_{t^\prime}^\prime$. 
Since $n_j^\prime$ and $w_2^\prime \cdots w_{t^\prime}^\prime$ are in $\{0, n-1\}^{<\omega}$, there must exist such an element $f^\prime$. 
By the construction, $g^\prime$ is in the subgroup that is isomorphic to $G_0$. 
In addition, each $(n-1)0n_j^\prime a_i$ is a fixed point of $g^\prime$, and $(n-1)0n_j^\prime b_i$ is in $\supp{g^\prime}$. 
By the construction of the map $I_{2, n}$ in Theorem \ref{prop_embedding_G_0(p)G_0(q)}, if we restrict the domain set of $g^\prime$ to $\{0, n-1\}^\omega$ and identify $n-1$ with $1$, then we obtain the element of $G_0$ which satisfies all the assumptions about fixed points and elements of support. 
However, since this does not happen by Proposition \ref{proposition_piecewiseprojective}, this is a contradiction. 

Finally, assume that $f(n_j)=w_1 \cdots w_t$ is in $\{0, n-1\}^{<\omega}$. 
Since $w_1y^{l_1}w_2 y^{l_2}\cdots w_t y^{l_t} a_i$ is in $\{0, n-1\}^{\omega}$, $n_j$ is also in $\{0, n-1\}^{<\omega}$. 
Therefore, $f$ is represented by a word on $\{x_0, \x{0}{(n-1)}\}$. 
We again note that some $l_q$ may be zero. 
Then $g^\prime=f\y{w_1\cdots w_t}{l_t} \cdots \y{w_1w_2}{l_2} \y{w_1}{l_1}$ is in $G_0(n)$, and an element of $G_0$ can be similarly constructed from $g^\prime$. 
\end{proof}
\end{lemma}
\begin{remark}
By the construction, the sequence $a_1, \dots, a_{2m}$ in Lemma \ref{lem_supp_finite_component} is uniquely determined. 
\end{remark}
The following lemma corresponds to \cite[Lemma 13]{monod2013groups}. 
\begin{lemma}
If $f, g \in G_0(n)$ have a common fixed point $x \in \NN$, then any element of the second derived subgroup $\langle f, g \rangle ^{\prime \prime}$ acts trivially on a neighborhood $\{\xi \mid \xi=a \xi^\prime, \xi^\prime \in \NN \}$ of $x$ where $a$ is in $\Nn$. 
\begin{proof}
Let $A$ be the set $\{a\xi \mid a \in \Nn, \xi \in \{0, n-1\}^\omega\}$. 
If $x \in A$, by appropriate conjugation, we regard $f$ and $g$ as elements in $G_0<G_0(n)$. 
Hence by \cite[Lemma 13]{monod2013groups}, the claim follows. 
If $x \not \in A$, we regard $f$ and $g$ as piecewise linear maps on $\mathbb{R}$. 
Then the first derived subgroup $\langle f, g \rangle ^{\prime}$ is trivial. 
\end{proof}
\end{lemma}
Using these lemmas, an argument similar to Monod's one in \cite[Theorem 14]{monod2013groups} shows the following: 
\begin{theorem}\label{Theorem_no_free_group}
The group $G_0(n)$ has no free subgroups. 
\end{theorem}
\begin{remark} \label{remark_torsion_free}
We can see that $G_0(n)$ is torsion free. 
Indeed, let $g$ $(g \neq \textrm{Id})$ be in $G_0(n)$ and assume that $g^m= \textrm{Id}$ holds. 
Since $g$ preserves the order of $\NN$, for $x$ in $\supp{g}$, we have either
\begin{align*}
&x < g(x)<g^2(x)< \cdots < g^m(x)=x
\shortintertext{or}
&x>g(x)>g^2(x)>\cdots >g^m(x)=x. 
\end{align*}
This is a contradiction. 
\end{remark}
\begin{remark}
If there exists an embedding from $G_0(n)$ into either $G_0$ or the Monod's group $H$ \cite{monod2013groups}, the theorem follows immediately because $G_0$ and $H$ contain no free subgroups. 
As mentioned in the introduction, we do not know whether $G_0(n)$ is a subgroup of $H$ or not. 
\end{remark}
\subsection{The abelianization of $G_0(n)$ and simplicity of the commutator subgroup}
The idea for the theorems in this section comes from \cite{brown1987finiteness} and \cite{burillo2018commutators}.

We use the surjective homeomorphism $a\colon F(n) \to \mathbb{Z}^n$ defined in \cite[section 4D]{brown1987finiteness}. 

\begin{theorem}[{\cite[Lemma 2.1]{burillo2018commutators}} for $n=2$] \label{Th_abelianization_G0(n)}
The abelianization of $G_0(n)$ is isomorphic to $\mathbb{Z}^{n+1}$. 
\begin{proof}
We define the map $\pi\colon Z(n) \to \mathbb{Z}^n \bigoplus \mathbb{Z}$ by
\begin{align*}
&\x{i}{s} \mapsto (a(\x{i}{s}), 0), & &y_s \mapsto (\bm{0}, 1). 
\end{align*}
Since $G_0(n)$ is an $(n+1)$-generated group, its abelianization is a quotient of $\mathbb{Z}^{n+1}$. 
So if we obtain a surjective homomorphism $G_0(n) \to \mathbb{Z}^{n+1}$, then this map must be its abelianization map. 
Thus it is sufficient to show that $\pi$ extends to a homomorphism $G_0(n) \to \mathbb{Z}^{n+1}$. 
In order to do this, we only need to check that the relations in $R(n)$ are satisfied, which is clear except for (4). 
By the definition of $a$, we have $a(\x{0}{s})=\bm{0}$ for each $\x{0}{s}$, where $s$ is in $\Nn$ such that $y_s$ is in $Y(n)$. 
Thus the relation (4) is also satisfied. 
\end{proof}
\end{theorem}
\begin{corollary}
Let $n, m \geq 2$. 
Then the groups $G_0(n)$ and $G_0(m)$ are isomorphic if and only if $n=m$ holds. 
\end{corollary}
In the following, we show that the commutator subgroup $G_0(n)^\prime=[G_0(n), G_0(n)]$ is simple. 
As in the case of $n=2$, first we show that the second derived subgroup $G_0(n)^{\prime \prime}=[G_0(n)^\prime, G_0(n)^\prime]$ is simple by using Higman's Theorem, and 
then we show that $G_0(n)^{\prime \prime}=G_0(n)^{\prime}$. 

Let $\Gamma$ be a group of bijections of a set $E$. 
For $\alpha \in \Gamma$, we write the set-theoretic support $\supp{\alpha}$ for the set $\{x \in E \mid \alpha(x) \neq x\}$. 
\begin{theorem}[{\cite[Theorem 1]{MR72136}}] \label{Theorem_Higman_simple}
Suppose that for every $\alpha, \beta, \gamma \in \Gamma \setminus \{ 1_\Gamma \}$, there exists $\rho$ such that 
\begin{align*}
\gamma \Bigl( \rho \bigl(\supp{\alpha}\cup \supp{\beta}) \Bigr) \cap \rho \bigl(\supp{\alpha}\cup \supp{\beta}\bigr) = \emptyset
\end{align*}
holds. 
Then the commutator subgroup $\Gamma^\prime$ is simple. 
\end{theorem}
\begin{theorem}[{\cite[Theorem 2]{burillo2018commutators}} for $n=2$] \label{theorem_commutator_simple}
The commutator subgroup of the group $G_0(n)$ is simple. 
\end{theorem}
The following two lemmas complete the proof. 
\begin{lemma}
The second derived subgroup $G_0(n)^{\prime \prime}$ is simple. 
\begin{proof}
Let $\alpha, \beta, \gamma \in G_0(n)^\prime$. 
Choose $x \in \supp{\gamma}$. 
If $\gamma(x)>x$, then let $I$ be the set $\{x^\prime \in \NN \mid x<x^\prime<\gamma(x)\}$. 
If $\gamma(x)<x$, then let $I$ be the set $\{x^\prime \in \NN \mid \gamma(x)<x^\prime<x\}$. 
Then we have $I \cap \gamma(I) = \emptyset$. 

We note that restrictions of $\alpha$ to $\{0^i \xi \mid \xi \in \NN\}$ and $\{(n-1)^i\xi \mid \xi \in \NN \}$ are in $F(n)^\prime$ for sufficiently large $i$. 
Hence by further replacing $i$ with sufficiently large $i^\prime$, restrictions of $\alpha$ to $\{0^{i^\prime} \xi \mid \xi \in \NN\}$ and $\{(n-1)^{i^\prime}\xi \mid \xi \in \NN \}$ are identity maps by the definition of the abelianization map $a$. 
The same holds for $\beta$. 
Then, by high transitivity of $F(n)^\prime$, there exists $\rho \in F(n)^\prime$ such that 
\begin{align*}
\rho \bigl(\supp{\alpha} \cup \supp{\beta}\bigr) \subset I
\end{align*}
holds. 
Therefore we have 
\begin{align*}
\rho \bigl(\supp{\alpha} \cup \supp{\beta} \bigr) \cap \gamma \Bigl(\rho \bigl(\supp{\alpha} \cup \supp{\beta} \bigr) \Bigr) \subset I \cap \gamma(I)=\emptyset. 
\end{align*}
By Theorem \ref{Theorem_Higman_simple}, $G_0(n)^{\prime \prime}$ is simple. 
\end{proof}
\end{lemma}
The next lemma is shown by replacing $1$ with $n-1$ in \cite[Proposition 2.5]{burillo2018commutators}. 
The argument can be simplified by assuming that $n>2$ holds. 
\begin{lemma}[{\cite[Proposition 2.5]{burillo2018commutators}} for $n=2$]
$G_0(n)^\prime=G_0(n)^{\prime \prime}$. 
\end{lemma}
\subsection{The center of $G_0(n)$}\label{subsec_center_G0(n)}
In this section, we show that the center of the group $G_0(n)$ is trivial. 
The idea for the theorem comes from \cite[Section 4]{cannon1996introductory}. 
Since we do not know whether $G_0(n)$ is a subgroup of Monod's group $H$, our proof is slightly different from \cite{burillo2018commutators}. 

Let $D(n)\coloneqq \{s \overline{0} \mid s \in \Nn\}$, where $\overline{0}$ denotes the element $000\cdots $ in $\NN$. 
Then the following holds. 
\begin{lemma} \label{lemma_rich_Thompson_F(n)}
For any $s\overline{0}$ in $D(n)$, there exists $x$ in $F(n)$ such that 
\begin{align*}
\supp{x}=(s\overline{0}, \overline{(n-1)})=\{\xi \in \NN \mid s\overline{0}<\xi< \overline{(n-1)}\}
\end{align*}
holds. 
\begin{proof}
If $s=0, \dots, (n-2)$, then $x_0, \dots, x_{n-2}$ satisfy the claim, respectively. 
By regarding $(n-1)\overline{0}$ as $(n-1)0\overline{0}$, we can assume that the length of $s$ is greater than or equal to $2$. 

Let $s=s^\prime i \overline{0}$ ($i \in \{0, \dots, n-2 \}$). 
If $s^\prime=(n-1)\cdots(n-1)$, then $\x{i}{s^\prime}$ satisfies the claim. 
If $s^\prime \neq (n-1)\cdots(n-1)$, by using $\x{i}{s^\prime}$, we can define an element in $F(n)$ as in Figure \ref{rich_elements_Fn}, which satisfies the claim. 
\begin{figure}[tbp]
	\centering
	\includegraphics[width=50mm]{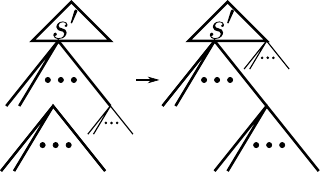}
	\caption{The construction of an element of $F(n)$ from $\x{i}{s^\prime}$. We add an $n$-caret to the leaf $s(n-1)$ of the domain tree and add an $n$-caret to the rightmost leaf of the range tree. }
	\label{rich_elements_Fn}
\end{figure}
\end{proof}
\end{lemma}
We note that $D(n)$ is a dense subset of $\NN$. 
\begin{theorem}[{\cite[Proposition 2.7]{burillo2018commutators}} for $n=2$] \label{theorem_center_trivial}
The center of $G_0(n)$ is trivial. 
\begin{proof}
Let $f$ be an element of the center of $G_0(n)$. 
For $g \in G_0(n)$, assume that $\supp{g}=(b_1, \overline{(n-1)})$ holds. 
Then we have $f(b_1)=b_1$. 
Indeed, if not, either $f(b_1)>b_1$ or $f^{-1}(b_1)>b_1$ holds. 
Note that $f(b_1)$ or $f^{-1}(b_1)$ is in $\supp{g}$. 
Since $g(b_1)=b_1$, in both cases, this contradicts that $fg=gf$ holds. 

By Lemma \ref{lemma_rich_Thompson_F(n)}, for every $s\overline{0} \in D(n)$, there exists $g \in F(n) \subset G_0(n)$ such that $b_1=s \overline{0}$ holds. 
Thus we have $f(s\overline{0})=s\overline{0}$ for every $s\overline{0} \in D(n)$. 
Since $D(n)$ is a dense subset of $\NN$, we conclude that $f$ is the identity map. 
\end{proof}
\end{theorem}
\subsection{Indecomposability with respect to direct products and free products}
In this section, we show that there exist no nontrivial ``decompositions'' using previous theorems.  

\begin{theorem}
There exists neither nontrivial direct product decompositions nor nontrivial free product decompositions of the group $G_0(n)$. 
\begin{proof}
Suppose that $G_0(n)$ is isomorphic to $K \times H$ for some groups $K$ and $H$. 
We first assume that $H$ is abelian. 
Then the center of $G_0(n)$ contains $\{1\} \times H$. 
Since the center of $G_0(n)$ is trivial (Theorem \ref{theorem_center_trivial}), $H$ must be the trivial group. 
From the same argument, if $K$ is abelian, then $K$ is trivial group. 

We assume that $K$ and $H$ are not abelian. 
We note that the commutator subgroup of $G_0(n)=K \times H$ is $[K, K] \times [H, H]$. 
Since $[K, K]$ and $[H, H]$ are not trivial, the group $\{1\} \times [H, H]$ is a nontrivial normal subgroup of $[K, K] \times [H, H]$. 
However, this contradicts that $[G_0(n), G_0(n)]$ is simple (Theorem \ref{theorem_commutator_simple}). 

Finally, we assume that $G_0(n)=K \star H$ for nontrivial groups $K$ and $H$. 
Let $k \in K \setminus \{1\}$ and $h \in H \setminus \{1\}$. 
Since $G_0(n)$ is torsion free (Remark \ref{remark_torsion_free}), both $k$ and $h$ generate infinite cyclic groups $\langle k \rangle$ and $\langle h \rangle$, respectively. 
This implies that $G_0(n)$ has a subgroup 
\begin{align*}
\langle k \rangle \star \langle h \rangle \cong \mathbb{Z} \star \mathbb{Z}=F_2. 
\end{align*}
By Theorem \ref{Theorem_no_free_group}, this is a contradiction. 
\end{proof}
\end{theorem}

\subsection*{Acknowledgments}
I would like to appreciate Professor Motoko Kato and Professor Shin-ichi Oguni for several comments and suggestions. 
I would also like to thank my supervisor, Professor Tomohiro Fukaya, for his comments and careful reading of the paper. 
Last but not least, I am grateful to the anonymous referees for providing many helpful comments and suggestions. 
\bibliographystyle{plain}
\bibliography{references} 

\address{
Department of Mathematical Sciences,
Tokyo Metropolitan University,
Minami-Osawa Hachioji, Tokyo, 192-0397, Japan
}

\textit{E-mail address}: \href{mailto:kodama-yuya@ed.tmu.ac.jp}{\texttt{kodama-yuya@ed.tmu.ac.jp}}
\end{document}